\documentclass{article}
\usepackage{row}
\usepackage{graphicx}
\begin{document}
\title{Key Developments in Geometry in the 19th Century
\footnote{This paper is dedicted to Mike Eastwood on the occaision of his sixtieth birthday in 2012 and to commemorate our collaboration in the early 1980s on the fascinating topic of twistor geometry.}
}
\author{Raymond O. Wells, Jr.
\footnote{Jacobs University, Bremen; University of Colorado at Boulder; wells@jacobs-university.de.}}

\maketitle
%\tableofcontents
\section{Introduction}
The notion of a manifold is a relatively recent one, but the theory of curves and surfaces in Euclidean 3-space \(\BR^3\) originated in the Greek mathematical culture. For instance, the book on the study of conic sections  by Apollonius  \cite{apollonius1896} described mathematically the way we understand conic sections today.  The intersection of a plane with a cone in \(\BR^3\) generated these curves, and Apollonius showed moreover that any plane intersecting  a {\it skew cone} gave one of the three classical conic sections (ellipse, hyperbola, parabola), excepting the degenerate cases of a point or intersecting straight lines. This was a difficult and important theorem at the time. The Greek geometers studied intersections of other surfaces as well, generating additonal curves useful for solving problems and they also introduced some of the first curves defined by transcendental functions (although that terminology was not used at the time), e.g., the quadratrix, which was used to solve the problem of squaring the circle (which they couldn't solve with straight edge and compass and which was shown many centuries later not to be possible).  See for instance Kline's very fine book on the history of mathematics for a discussion of these issues \cite{kline1972}.

Approximately 1000 years after the major works by the Greek geometers, Descartes published in 1637 \cite{descartes} a revolutionary book which contained a fundamentally new way to look at geometry, namely as the solutions of algebraic equations.  In particular the solutions of equations of second degree in two variables described precisely the conic sections that Apollonius had so carefully treated.  This new bridge between algebra and geometry became known during the 18th and 19th centuries as {\it analytic geometry} to distinguish it from {\it synthetic geometry}, which was the treatment of geometry as in the book of Euclid \cite{euclid}, the methods of which were standard in classrooms and in research treatises for approximately two millenia. Towards the end of the 19th century and up until today {\it algebraic geometry} refers to the relationship between algebra and geometry, as initiated by Descartes, as other forms of geometry had arisen to take their place in modern mathematics, e.g., topology, differential geometry, complex manifolds and spaces, and many other types of geometries.

Newton and Leibniz made their major discoveries concerning differential and integral calculus in the latter half of the 17th century (Newton's version was only published later; see \cite{newton1736}, a translation into English by John Colson from a Latin manuscript from 1671). Leibniz published his major work on calculus in 1684 \cite{leibniz1684}.  As is well known, there was more than a century of controversy about priority issues in this discovery of calculus  between the British and Continental scholars (see, e.g., \cite{kline1972}). In any event, the 18th century saw the growth of analysis as a major force in mathematics (differential equations, both ordinary and partial, calculus of variations, etc.). In addition the variety of curves and surfaces in \(\BR^3\) that could be represented by many new transcendental functions expanded greatly the families of curves and surfaces in \(\BR^3\) beyond those describable by solutions of algebraic equations. 

\begin{figure} 
\vspace{6pt}
\centerline{
	\includegraphics[width=12cm]{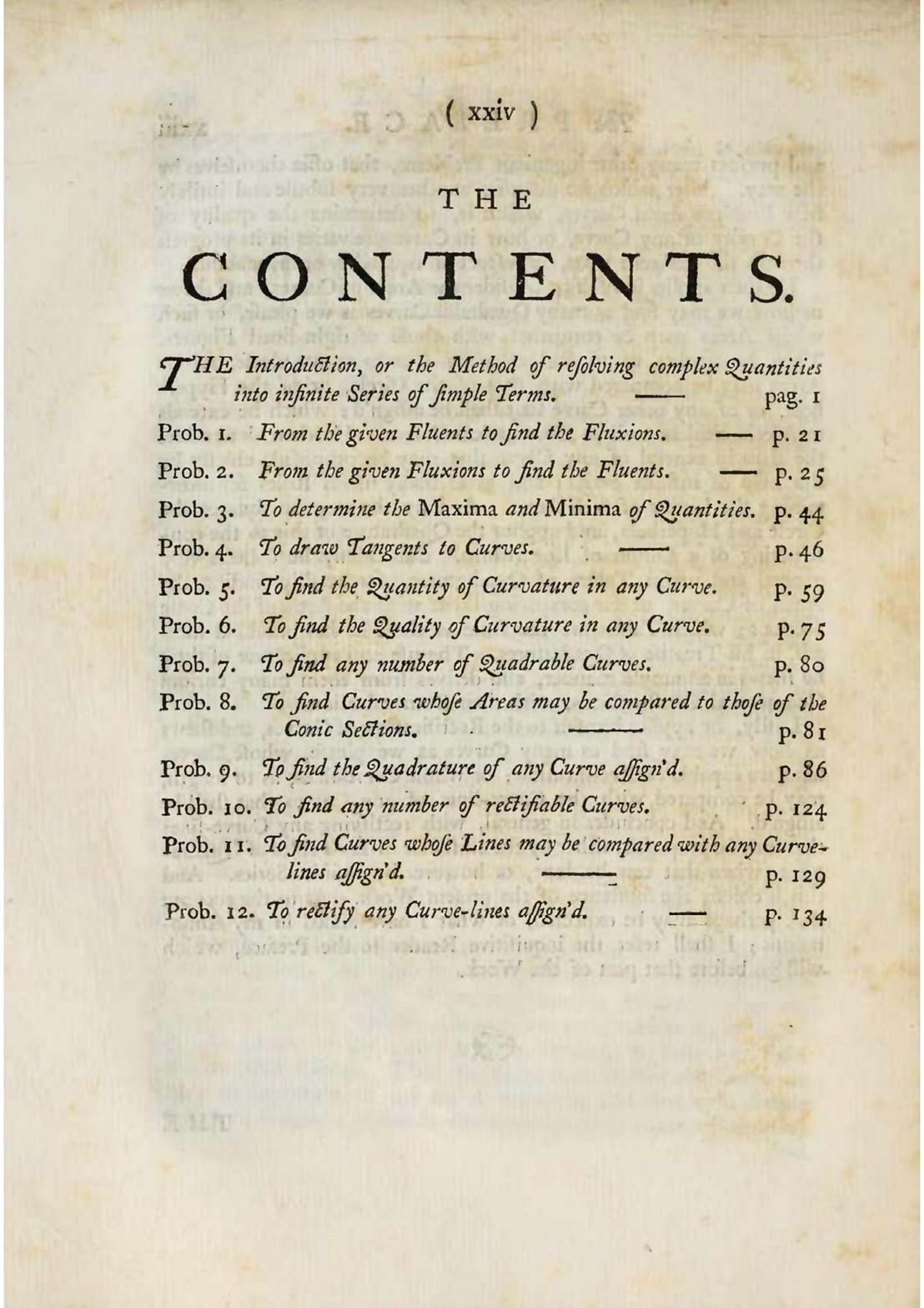}}
 \caption{Table of Contents of Newton's 1736 Monograph on Fluxions.}
 \label{fig:newton1736TOC}
\end{figure}

Moreover, a major development that grew up at that time, and which concerns us in this paper, is the growing interaction betwen analysis and geometry. An important first step was the analytic description of the curvature of a curve in the plane at a given point by Newton.  This was published in the 1736 monograph  \cite{newton1736}, which was a translation of his work from 1671.  We reproduce in Figure \ref{fig:newton1736TOC} the table of contents of this paper,
where we see that the study of curvature of a curve plays a central role.  In the text of this reference one finds the well-known formula for the curvature of a curve defined as the graph of a function \(y=f(x)\)
\be
\label{eqn:plane-curvature}
K_P=\pm \frac{f''(x)}{[1+(f'(x)^2]^\frac{3}{2}},
\ee
which one learns early on in the study of calculus.  Curvature had been studied earlier quite extensively by Huygens \cite{huygens1673}, and he computed the curvature of many explicit curves (including the conic sections, cycloids, etc.), but without calculus, however, using limiting processes of geometric approximations that are indeed the essence of calculus. Even earlier Apollonius had been able to compute the curvature of a conic section at a specific point (see \cite{apollonius1896} and Heath's discussion of this point).

The next development concerning analysis and geometry was the study of space curves in 
\(\BR^3\) with its notion of curvature and torsion as we understand it today.  This work on space curves was initiated by a very young (16 years old) Clairaut \cite{clairaut1731} in 1731.  Over the course of the next century there were numerous contributions to this subject by Euler, Cauchy, and many others, culminating in the definitive papers of Frenet and Serret in the mid-19th century \cite{serret1851}, \cite{frenet1852} giving us the formulas for space curves in \(\BR^3\) that we learn in textbooks today
\bea
\label{eq:curvature}
\BT'(s) & = &\k(s)\BN(s),\\
\label{eq:frenet-serret}
\BN'(s) & =&-\k(s)\BT(s) + \t(s)\BB,\\
\label{eq:torsion}
\BB'(s) & =& -\t(s)\BN(s),
\eea
where \(\BT(s)\), \(\BN(s)\) and \(\BB(s)\) denote the tangent, normal and binormal vectors to a given curve \(\Bx(s)\) parametrized by arc length \(s\). We need to recall that the vector notation we use today only came at the beginning of the 20th century, so the formulas of Serret and Frenet looked much more cumbersome at the time, but the mathematical essence was all there. Space curves were called "curves of double curvature" (French: {\it courbes \`{a} double courbure}) throughout the 18th and 19th centuries.

A major aspect of studying curves in two or three dimensions is the notion of the measurement of arc length, which has been studied since the time of Archimedes, where he gave the first substantive approximate value of \(\pi\) by approximating the arc length of an arc of a circle (see, e.g., \cite{resnikoff-wells1984}).  Using calculus one was able to formalize the length of a segment of a curve by the well-known formula
\[
\int_\G ds = \int^b_a \sqrt{(x'(t))^2+(y'(t))^2}dt,
\]
if \(\G\) is a curve in \(\BR^2\) parametrized by \((x(t),y(t)\) for \(a\le t \le b\). The expression \(ds^2= dx^2 + dy^2\) in \(\BR^2\) or \(ds^2 =dx^2 +dy^2 +dz^2\) in \(\BR^3\) became known as the {\it line element} or infinitesimal measurement of arc length for any parametrized curve.  What this meant is that a given curve could be aproximated by a set of secants, and measuring the lengths of the secants by the measurement of the length of a segment of a straight line in \(\BR^2\) or \(\BR^3\) (using Pythagoras!), and taking a limit, gave the required arc length. The line elements \(ds^2= dx^2 +dy^2\) and \(d^2 = dx^2 +dy^2 +dz^2\) expressed {\it both} the use of measuring straight-line distance in the ambient Euclidean space as well as the limiting process of calculus.  We will see later how generating such a line element to be {\it independent} of an ambient Euclidean space became one of the great creations of the 19th century.

\begin{figure} 
\vspace{6pt}
\centerline{
	\includegraphics[width=12cm]{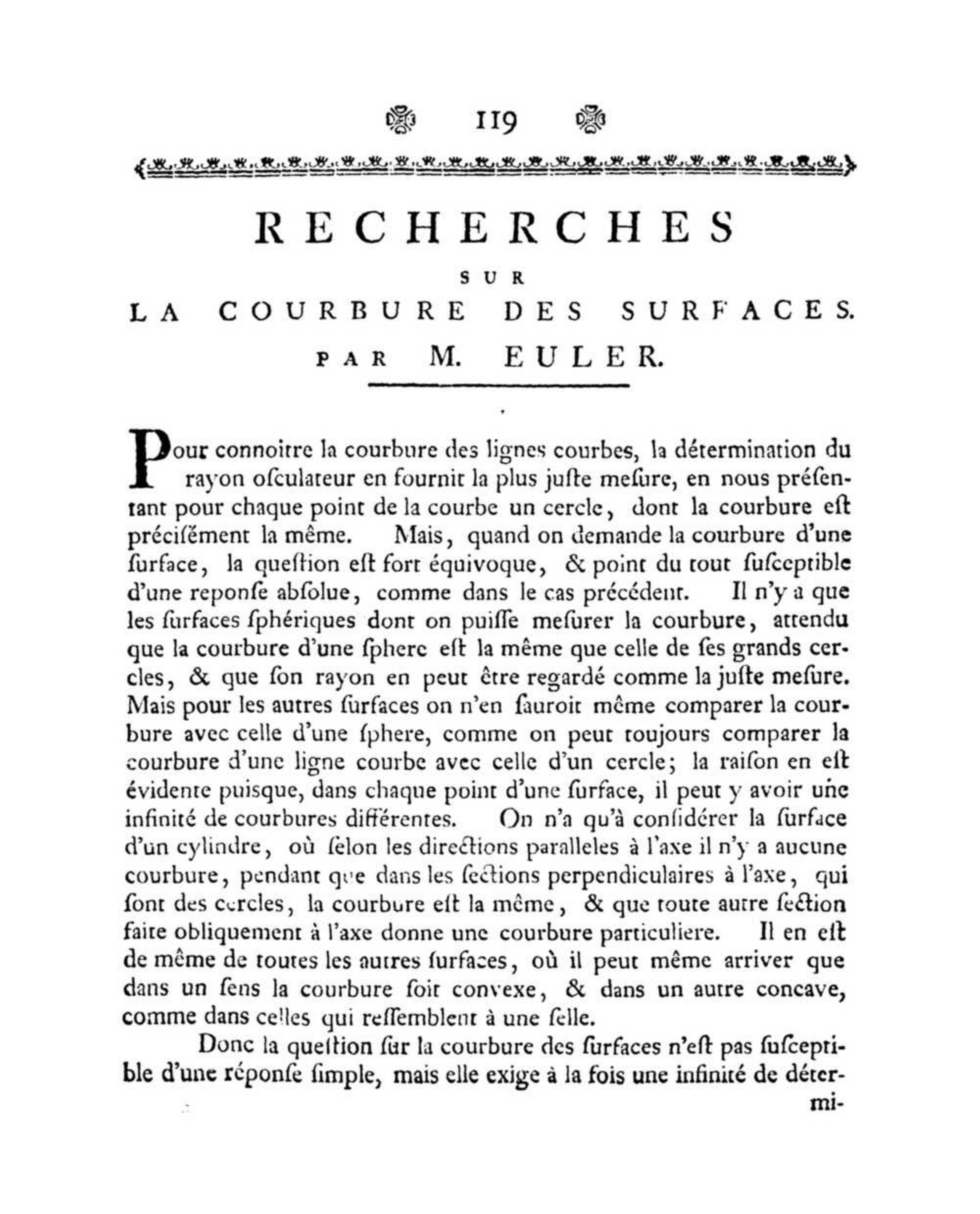}}
 \caption{The opening page of Euler's work on curvature.}
 \label{fig:euler1767p119}
\end{figure}

In our discussion above about curvature of a curve in \(\BR^2\) or \(\BR^3\) the basic understanding from Huygens to Frenet and Serret was that the curvature of a curve \(\G\) is the reciprocal of the radius of curvature of the osculating circle (best appoximating circle) in the osculating plane at a point on the parametrized curve \(\g(s)\) spanned by the tangent and normal vectors at \(\g(s)\).  The osculating circle to \(\G\) at \(\g(s)\) is a circle in the osculating plane that intersects \(\G\) tangentially at \(\g(s)\) and moreover also approximates the curve to second order at the same point (in a local coordinate Taylor expansion of the functions defining the curve and the circle, their Taylor expansions would agree to second order). This concept of curvature of a curve in \(\BR^3\) was well understood at the end of the 18th century, and the later work of Cauchy, Serret and Frenet completed this set of investigations begun by the young Clairaut a century earlier.  The problem arose:  how can one define the curvature of a surface defined either locally or globally in \(\BR^3\)?  This is a major topic in the next section, where we begin to survey more specifically the major developments in the 18th century by Gauss, Riemann and others concerning differential geometry on manifolds of two or more dimensions.  But first we close this introduction by quoting from Euler in his paper entitled ``Recherches sur la courbure des surfaces" \cite{euler1767} from 1767 (note this paper is written in French, not like his earlier works, all written in Latin). In Figure \ref{fig:euler1767p119} 
we show the first page of the article and we quote the translation here:
\begin{quote}
In order to know the curvature of a curve, the determination of the radius of the osculating circle furnishes us the best measure, where for each point of the curve we find a circle whose curvature is precisely the same. However, when one looks for the curvature of a surface, the question is very equivocal and not at all susceptible to an absolute response, as in the case above. There are only spherical surfaces where one would be able to measure the curvature, assuming the curvature of the sphere is the curvature of its great circles, and whose radius could be considered the appropriate measure. But for other surfaces one doesn't know even how to compare a surface with a sphere, as when one can always compare the curvature of a curve with that of a circle. The reason is evident, since at each point of a surface there are an infinite number of different curvatures. One has to only consider a cylinder, where along the directions parallel to the axis, there is no curvature, whereas in the directions perpendicular to the axis, which are circles, the curvatures are all the same, and all other obliques sections to the axis give a particular curvature. It's the same for all other surfaces, where it can happen that in one direction the curvature is convex, and in another it is concave, as in those resembling a saddle.
\end{quote}

In the quote above we see that Euler recognized the difficulties in defining curvature for a surfacs at any given point.  He does not resolve this issue in this paper, but he makes extensive calculations and several major contributions to the subject. He considers  a surface \(S\) in \(\BR^3\) defined as a graph
\[
z=f(x,y)
\]
near a given point \(P=(x_0,y_0,z_0)\).  At the point \(P\) he considers planes in \(\BR^3\) passing through the point \(P\) which intersect the surface in a curve in that given plane. For each such plane and corresponding curve he computes explicitly the curvature of the curve at the point \(P\) in terms of the given data.

He then restricts his attention to planes which are normal to the surface at \(P\) (planes containig the normal vector to the surface at \(P\)). There is a one-dimensional family of such planes \(E_\th\), parametrized by an angle \(\th\).  He computes explicitly the curvature of the intersections of \(E_\th\) with \(S\) as a function of \(\th\) , and observes that there is a maximum and minimum \(\k_\textrm{max}\) and \(\k_\textrm{min}\) of these curvatures at \(P\).%
\footnote
{Moreover, he shows that the plane \(E_{\th_\textrm{max}}\) is {\it perpendicular} to the plane \(E_{\th_\textrm{min}}\) with the minimum curvature \(\th_\textrm{min}\) (assuming the nondegenerate case).}
He also shows that if one knows the sectional curvature for angles \(\th_1\ne\th_2\ne\th_3\), then the curvature at any given angle is a computable linear combination of these three curvatures in terms of the given geometric data.  This is as far as he goes.  He does {\it not} use this data to define the {\it  curvature} of the surface \(S\) at the point \(P\).  This step is taken by Gauss in a visionary and extremely important paper some 60 years later, as we will see in the next section.

This is the background for our discussion in this paper of some key geometric discoveries in the 19th century, which completely transformed  the viewpoint of the 18th century of curves and surfaces in \(\BR^3\) (with the metric geometry induced from the Euclidean metric in \(\BR^3\)) to the contemporary 20th (and beginning of the 21st) century of abstract topological, differentiable and complex manifolds, with a variety of measuring tools on such manifolds.

What were the major geometric discoveries in the 19th century, and how did they come about?  The primary purpose of this historical paper is to present and discuss several of these major steps in mathematical creativity. The major topics we want to focus on include:
\bi
\item intrinsic differential geometry
\item projective geometry
\item higher-dimensional manifolds and differential geometry
\ei
What we won't try to cover in this paper is the development of complex geometry, including the theory of Riemann surfaces and projective algebraic geometry, as well as the discovery of set theory by Cantor and the development of topology as a new discipline. These ideas, the history of which is fascinating in its own right, as well as those topics mentioned above, led to the theory of topological manifolds and the corresponding theories of differentiable and complex manifolds, which became major topics of research in the 20th century. We plan to treat these topics and others in a book which is in the process of being developed.

\section{Gauss and Intrinsic Differential Geometry}
\label{sec:gauss}

The first major development at the beginning of the 19th century was the publication in 1828 by Gauss  of his historic landmark paper concerning differential geometry of surfaces entitled {\it Disquisitiones circa superficies curvas} \cite{gauss1828}.  He published the year before a very readable announcement and summary of his major results in \cite{gauss1827}, and we shall quote from this announcement paper somewhat later, letting Gauss tell us in his own words what he thinks the significance of his discoveries is.  For the moment, we will simply say that this paper laid the foundation for doing intrinsic differential geometry on a surface and was an important step in the creation of a theory of abstract manifolds which was developed a century later.

We now want to outline Gauss's creation of curvature of a surface as presented in his original paper. Consider a locally defined surface \(S\) in \(\BR^3\) which can be represented in either of three ways, as Gauss points out, and he uses all three methods in his extensive computations. We will follow his notation so that the interested reader could refer to the original paper for more details.  First we consider \(S\) as the zero set of a smooth function
\[
w(x,y,z)=0, w_x^2+w_y^2+w_z^2\ne0.
\]
Secoondly, we consider \(S\) as the graph of a function
\[
z=f(x,y),
\]
by the implicit function theorem.%
\footnote{Gauss used here simply \(z=z(x,y)\).}
And finally we consider \(S\) to be the image of a parmetric representation
\bean
x&=&x(p,q),\\
y&=&y(p,q),\\
z&=&z(p,q).
\eean
In fact, Gauss, as mentioned above, uses all three representations and freely goes back and forth in his computations, using what he needs and has developed previously. 

At each point \(P\) of the surface \(S\) there is a unit normal \(N_P\) associated with a given orientation of the surface, and we define the mapping
\[
g: S \rightarrow S^2 \subset \BR^3,
\]
where 
\(S^2\) is the standard unit sphere in \(\BR^3\) (\(S^2=\{(x,y,z): x^2+y^2+z^2=1\}\)), and where \(g(P):= N_P\).  This mapping, first used in this paper, is now referred to as the {\it Gauss map}. If \(U\) is an open set in \(S\), then the area (more precisely measure) on the 2-sphere of \(g(U)\) is called by Gauss the {\it total curvature} of \(U\).  For instance if \(U\) is an open set in a plane, then all the unit normals to \(U\) would be the same, \(g(U)\) would be a single point, its measure would be zero, and thus the total curvature of the planar set would be zero, as it should be.

Gauss then defines the {\it curvature} of \(S\) at a point \(P\) to be the derivative of the mapping \(g:S\rightarrow S^2\) at the point \(P\).  He means by this the ratio of infinitesimal areas of the image to the infinitesimal area of the domain, and this is, as Gauss shows, the same as the Jacobian determinant of the mapping in classical calculus language, or in more modern terms the determinant of the derivative mapping \(dg: T_P(S)\rightarrow T_{g(P)}S^2\). Today we refer to the curvature at a point of a surface as defined by Gauss to be {\it Gaussian curvature}.  Note that this definition is {\it a priori} extrinsic in nature, i.e., it depends on the surface being embedded in \(\BR^3\) so that the notion of a normal vector to the surface at a point makes sense. Gauss proceeds to compute the curvature of a given surface in each of the three representations above. In each case he expresses the normal vector in terms of the given data and explicitly computes the required Jacobian determinant. We will look at the first and third cases in somewhat more detail.

We start with the simplest case of a graph
\[
z=f(x,y)
\]
and let 
\[
df=tdx+udy,  
\]
where \(t=f_x, u=f_y.\)
Similarly, let \(T :f_{xx}, U := f_{xy}, V:= f_{yy}\), and we have
\bean
dt&=&Tdx+Vdy,\\
du&=&Udx+Vdy.
\eean
Thus the unit  normal vectors to \(S\) have the form
\[
(X,Y,Z)=(-tZ, -uZ,Z),
\]
where 
\[
Z^2 = \frac{1}{1+t^2+u^2}.
\]
By definition the curvature is the Jacobian determinant
\[
k=\frac{\partial X}{\partial x}\frac{\partial x}{\partial y}- \frac{\partial Y}{\partial x}\frac{\partial Y}{\partial y},
\]
which Gauss computes to be 
\[
k=\frac{TV-U^2}{(1+t^2+u^2)},
\]
or in terms of \(f\),
\be
\label{eqn:curvature-graph}
k=\frac{f_{xx}f_{yy}-f_{xy}^2}{(1+f_x^2+f_y^2)^2},
\ee
and we see a strong similarity to the formula for the curvature of a curve in the case of a graph of a function as given by Newton (\ref{eqn:plane-curvature}).

Gauss considers the case where the tangent plane at \(x_0,y_0\) is perpendicular to the \(z\)-axis (i.e., \(f_x(x_0,y_0)=f_y(x_0,y_0)= 0\)), obtaining 
\[
k=f_{xx}f_{yy}- f_{xy}^2,
\]
and by making a 
rotation in the \((x,y)\) plane to get rid of the term \(f_{xy}(x_0,y_0)\), he obtains
\[
k=f_{xx}f_{yy},
\]
which he identifies as being the product of the two principal curvatures of Euler.

Next Gauss proceeds to compute the curvature of S in terms of an implicit representation of the form 
\(w(x,y,z)=0\), and he obtains a complicated expression which we won't reproduce here, but it has the form
\[
(w_x^2+w_y^2+ w_z^2)k=h(w_x,w_y,w_z,w_{xx},w_{xy},w_{yy},w_{xz},w_{yz},w_{zz}),
\]
where \(h\) is an explicit homogeneous polynomial of degree 4 in 9 variables.  Here, of course, \(w_x^2+w_y^2+w_z^2\ne 0\) on \(S\), so this expresses \(k\) as a rational function of the derivatives of first and second order of \(w\), similar to (\ref{eqn:curavture-graph}) above. This explicit formula is given on p. 232 of \cite{gauss1828}.

Then Gauss moves on to the representation of curvature in terms of a parametric representation of the surface.  He first computes the curvature in an extrinsic manner, just as in the preceding two cases, and then he obtains a new representation, by a very slick change of variables, which becomes the intrinsic form of Gaussian curvature. Here are his computations.  He starts with the representation of the local surface as 
\bean
x&=&x(p,q),\\
y&=&y(p,q),\\
z&=&z(p,q),
\eean
and he gives notation for the first and second derivatives of these functions.  Namely, let 
\bean
dx&=&adp+a'dq,\\
dy&=&bdp+b'dq,\\
dz&=&cdp+c'dq,
\eean
where, of course, \(a= x_p(p.q), b =y_p(p,q)\), etc., using the subscript notation for partial derivatives. In the same fashion we define
the second derivatives
\[
\ba{lll}
\a:=x_{pp}, & \a':=x_{pq}, & \a'':= x_{qq},\\
\b:= y_{pp}, & \b' := y_{pq}, & \b'':= y_{qq},\\
\g:=z_{pp}, & \g':= z_{pq}, & \g'' := z_{qq}.
\ea
\]
Now the vectors \((a,b,c), (a',b', c')\) in \(\BR^3\) represent tangent vectors to \(S\) at \((x(p,q), y(p,q), q(p,q))\in S\), and thus the vector
\be
\label{eqn:normal}
(A,B,C):= (bc'-cb',ca'-ac',ab'-ba')
\ee
is a normal vector to S at the same point (cross product of the two tangent vectors).  Hence one obtains
\[
Adx+Bdy + Cdz = 0
\]
on \(S\), and therefore we can write 
\[
dz= -\frac{A}{C}dx -\frac{B}{C}dy,
\]
assuming that \(C\ne0\), i.e., we have the graphical representation \(z= f(x,y)\), as before.
Thus
\bea
\label{eqn:t}
z_x = t = -\frac{A}{C},\\
\label{eqn:u}
z_y= u= -\frac{B}{C}.
\eea
Now since 
\bean
dx&=&adp+a'dq,\\
dy&=& bdb +b'dq,
\eean
we can solve these linear equations for \(dp\) and \(dq\), obtaining (recalling that \(C\) is defined in (\ref{eqn:normal}))
\bea
\label{eqn:dp}
Cdp&=&b'dx-a'dy,\\
\label{eqn:dq}
Cdq&=& -bdx+ady.
\eea
Now we differentiate the two equations (\ref{eqn:t}) and (\ref{eqn:u}), and using (\ref{eqn:dp}) and (\ref{eqn:dq}) one obtains
\[
C^3dt= (A\frac{\partial C}{\partial p}-C\frac{\partial A}{\partial p})(b'dx-a'dy) + (C\frac{\partial A}{\partial q}-A\frac{\partial C}{\partial q})(bdx-ady),\]
and one can derive a similar expression for \(C^3du\).

Using the expression for the curvature in the graphical case (\ref{eqn:curvature-graph}), and by setting
\bean
D&=& A\a + B\b +C\g,\\
D'&=& A\a' + B\b' + C\g',\\
D'' &=& A\a'' +B\b'' +C\g'',
\eean
Gauss obtains the formula for the curvature in this case
\be
\label{eqn:curvature-parameter}
k=\frac{DD'' - (D')^2}{(A^2+B^2+ C^2)^2},
\ee
where again this is a rational function of the first and second derivatives of the parametric functions \(x(p,q), y(p,q), z(p,q)\). Note that the numerator is homogenous of degree 4 in this case, and again this formula depends on the explicit representation of \(S\) in \(\BR^3\). One might think that he was finished at this point, but he goes on to make one more very ingenious  change of variables that leads to his celebrated discovery.

He defines
\bean
E&=& a^2 +b^2 + c^2,\\
F&=& aa' + bb' + c c',\\
G&=& (a')^2 + (b')^2 +(c')^2\\
\D&=& A^2 + B^2 +C^2=EG-F^2,
\eean
and additionally,
\bean
m&=& a\a +b\b +c\g,\\
m'&=& a\a' +b\b' + c\g',\\
m''&=& a\a'' +b\b'' + c \g'',\\
n&=& a'\a + b'\b + c'\g,\\
n'&=& a'\a' + b'\b' + c'\g',\\
n''&=& a'\a'' + b'\b'' +c'\g''.
\eean
He is able to show that (using the subscript notation for differentiation again, i.e., \(E_p=\frac{\partial E}{\partial p}\), etc.) 
\[
\ba{rclrcl}
E_p&=&2m,&\;E_q&=&2m',\\
F_p&=&m'+n,&\; F_q&=& m'+n',\\
G_p&=&2n', &\; G_q&=&2n'',
\ea
\]
from which one obtains, by solving these linear equations
\[
\ba{rclrcl}
m&=&\frac{1}{2}E_p, &\;n&=& F_o-\frac{1}{2}E_q,\\
m'&=&\frac{1}{2}E_q,& \;n'&=& \frac{1}{2}G_p,\\
m''&=&F_q-\frac{1}{2}G_p,&\; n'' &=&\frac{1}{2}G_q.
\ea
\]
Gauss had also shown that the numerator which appears in (\ref{eqn:curvature-parameter}) can be expressed in terms of these new variables as
\bea
DD'' -(D')^2 &=&\D[\a\a''-\b\b''+\g\g''- (\a')^2-(\b')^2-(\g')^2]\nonumber\\
&&+E([(n')^2-nn''] +F[nm'' -2m'n' +mn'']\nonumber\\
&&+G[(m')^2-mm''].\label{eqn:numerator}
\eea
Finally one confirms that 
\bea
\frac{\partial n}{\partial q}-\frac{\partial n'}{\partial p} &=& \frac{\partial m''}{\partial p}- \frac{m'}{\partial q}\nonumber\\
&=& \a\a''-\b\b'' + \g\g''-(\a')^2-(\b')^2-(\g')^2. \label{eqn:ii}
\eea
Substituting  (\ref{eqn:ii}) and (\ref{eqn:numerator}) into the curvature formula (\ref{eqn:curvature-parameter}) Gauss obtains the expression for the curvature he was looking for:
\bea 4(EG-FF)^2 k &=& E(E_qG_q -2F_pG_q + G_p^2)\nonumber\\
&&+F(E_pG_q-E_qG_p-2E_qF_q+4F_pF_q-2F_pG_p)\nonumber\\
&& +G(E_pG_p-2E_pF_q+E_q^2)\nonumber\\
&& -2(EG-FF)(E_{qq}-2F_{pq} +G_{pp}.\label{eqn:egregrium}
\eea

If we consider the metric on \(S\) induced by the metric on \(\BR^3\), we have,
in terms of the paramtric representation,
\bea
ds^2 &=& dx^2 +dy^2+ dz^2\nonumber\\
&=& (adp+a'dq)^2 +(bdp+b'dq)^2 +(cdp+c'd1)^2\nonumber\\
&=& Edp^2 +2Fdpdq+Gdq^2.\label{eqn:metric}
\eea
We can now observe that in (\ref{eqn:egregrium}) Gauss has managed to express the curvature of the surface in terms of \(E, F, G\) and its first and second derivatives with respect to  the parameters of the surface.  That is, the curvature is a function of the {\it line element} (\ref{eqn:metric}) and its first and second derivatives. 

\begin{figure} 
\vspace{6pt}
\centerline{
	\includegraphics[width=12cm]{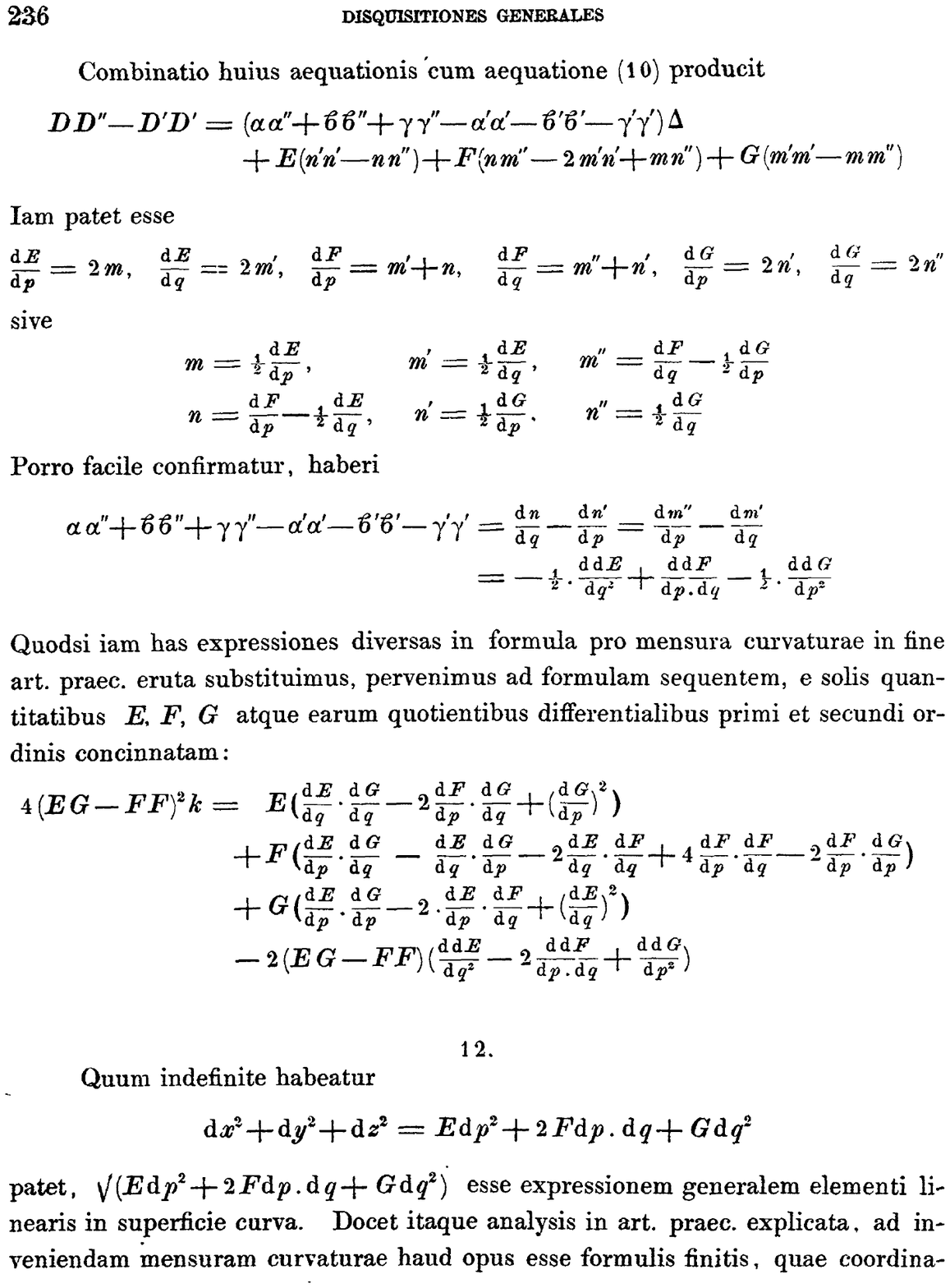}}
 \caption{The formula for Gaussian curvature in its intrinsic form.}
 \label{fig:curvature}
\end{figure}

\begin{figure} 
\vspace{6pt}
\centerline{
	\includegraphics[width=12cm]{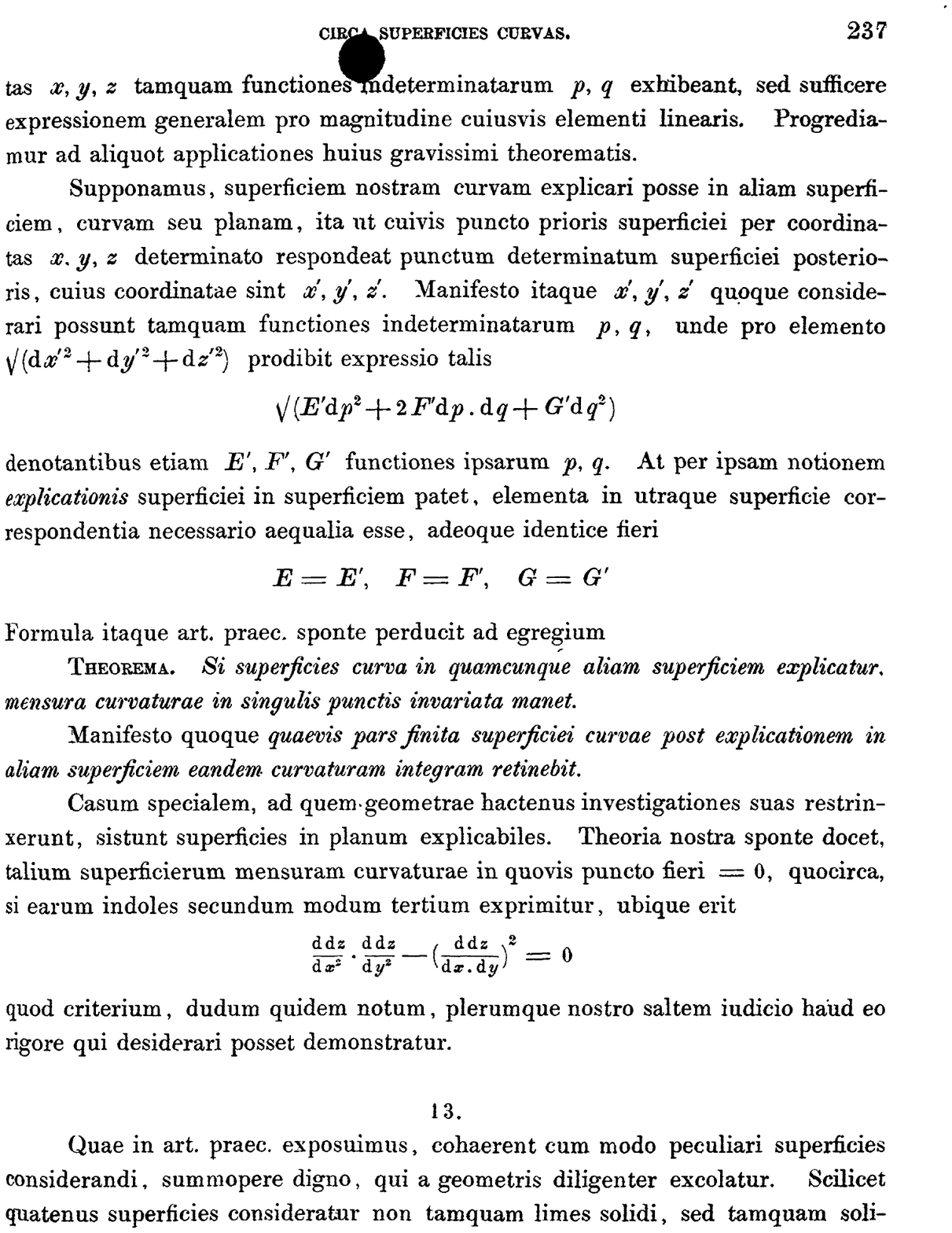}}
 \caption{The {\it Theorema Egregrium} of Gauss.}
 \label{fig:egregrium}
\end{figure}

Gauss called this his {\it Theorema Egregrium} (see pages 236 and 237 from \cite{gauss1828} given in Figures \ref{fig:curvature} and \ref{fig:egregrium}, where the formula is on p. 236 and the statement of the theorem is on page 237).  
What the text on p. 237 of Gauss's paper says, in so many words, is that if two surfaces can be represented by the same parametrization such that the induced metrics are the same, then the curvature is preserved.  In more modern language, an isometric mapping of one surface to another will preserve the Gaussian curvature. 

Gauss understood full well the significance of his work and the fact that this was the beginning of the study of a new type of geometry (which later generations have called {\it intrinisic differential geometry}).  We quote here from his announcement of his results published some months earlier (in his native German) from pp. 344--345 of \cite{gauss1827}:
\begin{quote}
\label{gauss-quote}
Diese S\"{a}tze f\"{u}hren dahin, die Theorie der krummen Fl\"{a}chen aus einem neuen Gesichtspunkte zu betrachten, wo sich der Untersuchung ein weites noch ganz  unangebautes Feld \"{o}ffnet.  Wenn man die Fl\"{a}chen nicht als Grenzen von K\"{o}rpern, sondern als K\"{o}rper, deren eine Dimension verschwindet, und zugleich als biegsam, aber nicht als dehnbar betrachte, so begreift man, dass zweierlei wesentlich verschiedene Relationen zu unterscheiden sind, theils nemlich solche, die eine bestimmte Form der Fl\"{a}che im Raum voraussetzen, theils solche, welche von den verschiedenen Formen, die die Fl\"{a}che annehmen kann, unabh\"{a}ngig sind. Die letz\-tern sind es, wovon hier die Rede ist: nach dem, was vorhin bemerkt ist, geh\"{o}rt dazu das Kr\"{u}mmungsmaass; man sieht aber leicht, dass eben dahin die Betrachtung der auf der Fl\"{a}che construirten Figuren,  ihrer Winkel, ihres Fl\"{a}cheninhalts und ihrer Totalkr\"{u}mmung, die Verbindung der Punkte durch k\"{u}rzesten Linien u. dgl. Geh\"{o}rt.\footnote{These theorems lead us to consider the theory of curved surfaces from a new point of view, whereby the investigations open to a quite new undeveloped field.  If one doesn't consider the surfaces as boundaries of domains, but as domains with one vanishing dimension, and at the same time as bendable but not as stretchable, then one understands that one needs to differentiate between two different types of relations, namely, those which assume the surface has a particular form in space, and those that are independent of the different forms a surface might take.  It is this latter type that we are talking about here.  From what was remarked earlier, the curvature belongs to this type of concept, and moreover figures constructed on the surface, their angles, their surface area, their total curvature as well as the connecting of points by curves of shortest length, and similar concepts, all belong to this class.}
\end{quote}
The results described in this short announcement (and the details in the much longer Latin papers on the subject) formed the basis of most of what became modern differential geometry.  An important point that we should make here is that Gauss did significant experimental work on measuring the curvature of the earth in the area around G\"{o}ttingen, where he spent his whole scientific career.  This involved measurements over hundreds of miles, and involved communicating between signal towers from one point to another.  He developed his theory of differential geometry as he was conducting the experiments, and at the end of all of his papers on this subject one finds descriptions of the experimental results (including, in particular, the short paper in German \cite{gauss1827}).

\section{Projective Geometry}
\label{sec:projective geometry}
The second major development that had its beginnings at the start of the 19th century was the creation of projective geometry which led in the latter half of the 19th century to the important concept of {\it projective space} \(\BP_n\) over the real or complex numbers (or more general fields).  This became a centerstage for the developments in algebraic geometry and complex manifolds at that time which carried over vigorously into the 20th century. In this section we will outline the very interesting story of how and why projective geometry developed.  

Projective geometry started as a school of mathematics in France around 1800.  Throughout the 18th century, as we mentioned above, Cartesian geometry and its interaction with differential and integral calculus, what became known later as differential geometry, dominated mathematical research in geometry.  In particular, the classical ideas of what became known as {\it synthetic geometry} in the spirit of Euclid's {\it Elements} began to fall to the wayside in mathematical research, even though Newton and others had often resorted to synthetic geometric arguments early in the 18th century as a complement to (and often as a check on) the analytic methods using coordinates. Some of the proponents of projective geometry wanted to create a type of geometry that could play an important alternative to coordinate geometry for a variety of interesting geometric problems.

It is generally recognized that projective geometry is the singular creation of Gaspard Monge
\footnote{Monge was a major scientific advisor to Napoleon on his Egyptian expedition and helped create the Ecole Polytechnique, the first engineering school in France, as well as the metric system.} and his pupils, principally Carnot, Poncelet,  and Chasles, to mention just three of them.  Monge published a book \cite{monge1827} {\it Geom\'{e}trie Descriptive} which inspired major developments by his pupils which, in turn, resulted in several substantive books, Poncelet's book in 1822 {\it Propriet\'{e}s Projectives}  \cite{poncelet1865}  being perhaps the most influential. The book by Chasles \cite{chasles1837} is a two-part work, the first of which is a brilliant historical treatise on the whole history of geometry up to that point in time, 1831. The second part is his own treatment of projective geometry, following up on the work of Monge, Carnot, Poncelet and others. Carnot's book {\it G\'{e}om\'{e}trie de Position} \cite{carnot1803}, along with Monge's original book from 1799 \cite{monge1827}, was a major influence on Poncelet.  We will  say more about this later, but first we want to give a overview of how and why projective geometry came to be and what these first authors believed they had achieved.

It is quite fascinating to read the 1827  edition of Monge's book \cite{monge1827}, which was edited by his pupil Bernab\'{e} Brisson and published after Monge's death in 1818. This is the fifth edition of the book which first appeared in 1799 and was the result of his lectures at the {\it \'{E}cole Normale},which were  very inspirational, according to several testimonials published in this edition. Monge was very concerned about secondary education and the creation of a new generation of educated citizens who could help in the development of the Industrial Revolution in France. He believed that descriptive geometry, which he developed in this book, would be a tool for representing three-dimensional objects in terms of their projections onto one or more planes in three-dimensional space, and that this should be a major part of the educational development of students. The applications he uses as examples came from architecture, painting, the representation of military fortifications, and many others. Monge was also interested in developing an alternative approach to classical synthetic geometry in the three-dimensional setting in order to develop a geometric method which would be useful in engineering.. Today we use his ideas in the form of blueprints for industrial design with such two-dimensional drawings representing  horizontal and vertical projections of the object being designed or manufactured.

The basic thesis of plane geometry as formulated by the Greeks and brought down to us in the book of Euclid is a solution of geometric problems in the plane by the use of the straight edge and compass. Monge's fundamental thesis in his book is to reduce problems of geometry in three-dimensional space to plane geometry problems (in the Euclidean sense) on the {\it projections} of the problems to two (or more) independent planes. Considering the straight edge and compass as engineering tools, a designer or engineer could work on three-dimensional problems in the two-dimensional medium.

Let us illustrate this with a simple example from his book.
\begin{figure}
\vspace{6pt}
\centerline{
	\includegraphics[width=8cm]{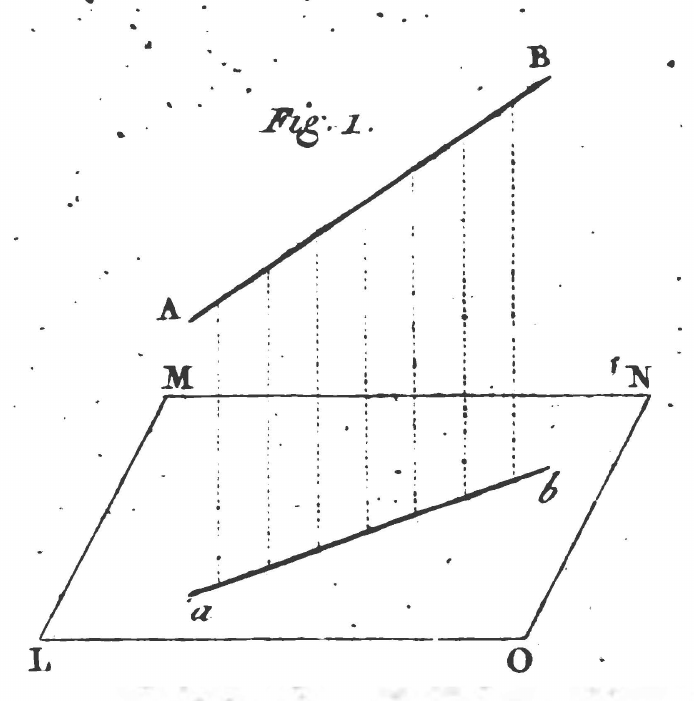}}
 \caption{Figure 1 in Monge's {\it G\'{e}om\'{e}trie Descriptive} \cite{monge1827}.}
 \label{fig:monge-fig1}
\end{figure}
In Figure \ref{fig:monge-fig1} we see the projection of the line segment \(AB\) in \(\BR^3\) onto the line segment \(ab\) in the plane \(LMNO\). Let us visualize these line segments as representing straight lines extended infinitely in both directions. Let's change notation slightly from that used by Monge in this figure. Let \(L\) be the given line in \(\BR^3\)  and let \(E_1\) and \(E_2\) be two non-parallel planes in \(\BR^3\). Let \(L_1\) and \(L_2\) be the perpendicular projections of the line \(L\) onto the planes \(E_1\) and \(E_2\), respectively. We see that \(L_1\) and \(L_2\) represent \(L\) in this manner. Can we reverse the process? Indeed, given \(L_1\) and \(L_2\), two lines in \(E_1\) and \(E_2\), and let \(H_1\) and \(H_2\) be the planes in \(\BR^3\) which are perpendicular to \(E_1\) and \(E_2\) and which pass through the lines \(L_1\) and \(L_2\). Then the intersection \(H_1\cap H_2\) is the desired line in \(\BR^3\) whose projections are \(L_1\) and \(L_2\).  This simple example is the main idea in the projection of a point in \(\BR^3\) onto the three orthogonal axes (\(x\)-axis,  \(y\)-axis, and \(z\)-axis), giving the Cartesian coordidnates \((x,y,z)\) of a given point. 

Monge formulates and solves a number of problems using these types of ideas.  Here is a quite simple example from his book, paralleling a similar question in plane geometry.  Given a line \(L\) in \(\BR^3\) and a point \(P\) not on the line, construct a line through \(P\) parallel to \(L\).  One simply chooses two reference planes \(E_1\) and \(E_2\), as above, and projects orthogonally \(L\) and \(P\) onto the planes, obtaining lines \(L_1\) and  \(L_2\) and points \(P_1\), and \(P_2\) in the reference planes.  Then, in this plane geometry setting, find parallel lines \(l_1\subset E_1\) and \(l_2\subset E_2\) to the lines \(L_1\) and \(L_2\) passing through the points \(P_1\) and \(P_2\).  The lines \(l_1\) and \(l_2\) determine a line \(l \subset \BR^3\), which solves the problem.

Monge considers in great depth for most of his book much more complicated problems concerning various kinds of surfaces in \(\BR^3\), and shows how to construct tangents, normals, principal curvatures and solutions to many other such problems. He uses consistently the basic idea of translating a three-dimensional problem into several two-dimensional problems.  This is the essence of descriptive geometry, but it becomes the basis for the later developments in projective geometry, as we will now see.

The most influential figure in the development of projective geometry in the first half of the 19th century (aside from the initial great influence of Monge, as described above) was undoubtedly Poncelet . He was an engineer and mathematician who served as Commandant of the Ecole polytechnique in his later years. His major contributions to projective geometry were his definitive books {\it Trait\'{e} des Propriet\'{e}s Projectives des Figures}, Volumes 1 and 2 \cite{poncelet1865}, \cite{poncelet1866}, which were first published in 1822 and 1824 and reappeared as second editions in 1865 and 1866.  In addition he published in 1865 a fascinating historical book {\it Applications d'Analyse et de G\'{e}om\'{e}trie} \cite{poncelet1862}, the major portion of which is the reproduction of notebooks Poncelet wrote in a Russian prisoner-of-war camp in Saratov 1813--1814.  He was interned there for about two years after a major battle (at Krasnoi) which Napoleon lost in November of 1812 towards the end of his disastrous Russian campaign to Moscow, and where Poncelet and others had been left on the field as dead. Poncelet had been serving
 as an {\it Officier de genie} (engineering officer). He had no books or notes with him and wrote out and developed further ideas he had learned from the lectures and writings of Monge and Carnot, in particular the works cited above \cite{monge1827} and \cite{carnot1803}. He attributes the long time between his first book (1822) and their revisions and his historical book which appeared some 40 to 50 years later to the extensive time commitment his administrative career demanded.  He also published numerous mathematical papers early in his career, most of which play an important role in his books on geometry. Moreover he was the author of several engineering monographs as well.

Poncelet promoted in his writings a specific doctrine for the development of geometry which became, with time, known as {\it projective geometry}, a term we didn't see him using in our reading of his works. He always used the expression ``{\it propriet\'{e}s projectives}" (projective properties) of figures, by which he meant those geometric properties of geometric figures that could be derived via projection methods by the new methodology that he was developing. He believed that the use of coordinate systems for the study of geometric problems was overvalued and led to problems of geometric understanding when negative and imaginary (complex) numbers appeared as solutions to equations. A real positive number could represent the length of a segment, or an area or volume of the figure, but what did negative and complex numbers represent geometrically? Only in the latter half of the 19th century did satisfactory answers to these questions arise.

As an example of this questioning of the use of algebra and geometry it is interesting to point to the opening 30 pages or so of the book by Carnot \cite{carnot1803}, which is solely dedicated to showing that negative numbers do not exist. Towards the end of his polemeical assertions about the non-existence of negative numbers, Carnot introduces the notion of ``direct" and ``inverse" of numbers to justify the proliferation of plus and minus signs (the usual formulas of algebra and trigonometry) in his book.  He also gives tantalizing hints of what became ``analysis situs" in the late 19th century, which is now called simply {\it topology}. This early work of Carnot included specifically the notion of points and lines at infinity as well as the notion of duality and dual problems, which we will say more about later.

At the time of Apollonius, we saw that any section of a skew cone could be considered as a section of a right circular cone, i.e. one of the classical conic sections. Now let's use the same picture but from a different perspective (to use a pun, in this case intended). Let  \(E_1\) and \(E_2\) be two (non-parallel) planes in \(\BF^3\) and let \(\G_1\) be a circle in \(E_1\) and let \(P\) be a point not on either plane, and let \(C\) be the cone formed by the pencil of lines emanating from \(P\), the {\it point of perspective}, and passing through the points of \(\G_2\).  Then consider the intersection \(\G_2\) of the cone \(C\) with the second plane \(E_2\). Then, according to Apollonius, this curve is again a conic section.  From the point of view of projective geometry one says that the curves \(\G_1\) and \(\G_2\) are projectively equivalent, and that they represent the ``same curve" projected onto different planes.  One can imagine a figure being projected onto one plane from one perspective point and also being in the same plane from a different perspective point, and then from the second perspective point being projected onto a third version of the figure on a different plane.  All of these changes of perspective and projections correspond to the natural mappings of projective space onto itself in modern terms, but this was the initial way these things were looked at by Poncelet and his school. 

\begin{figure}
\vspace{6pt}
\centerline{
	\includegraphics[width=10cm]{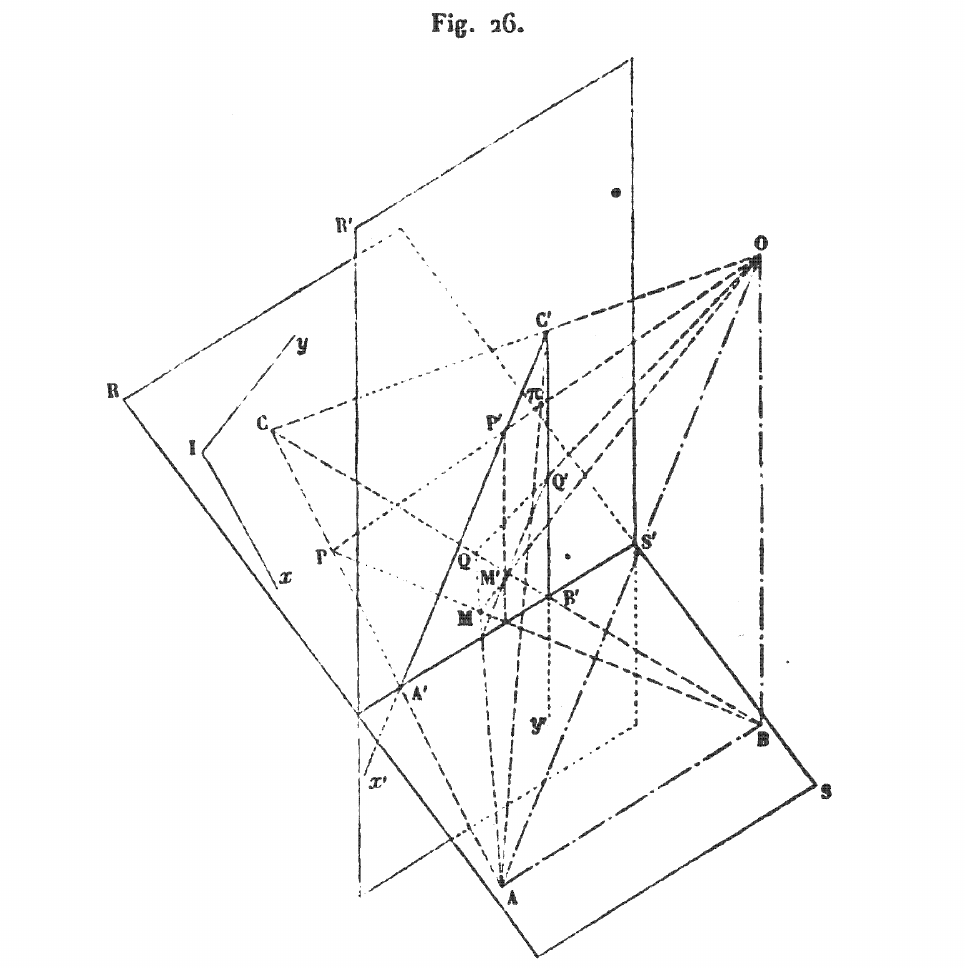}}
 \caption{Figure 26 in Poncelet's {\it Applications d'Analyse et G\'{e}om\'{e}trie.} \cite{poncelet1862}.}
 \label{fig:poncelet1862-p512}
\end{figure}

\begin{figure}
\vspace{6pt}
\centerline{
	\includegraphics[width=12cm]{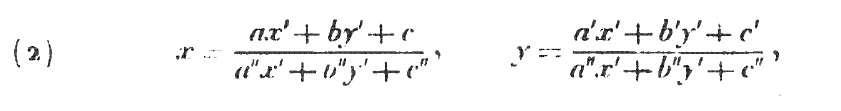}}
 \caption{Equation 2 on p. 616 in Poncelet's Figure 1 in {\it Applications d'Analyse et G\'{e}om\'{e}trie} \cite{poncelet1862}.}
 \label{fig:poncelet1862-p516}
\end{figure}

They took the time to represent these changes of perspective and projection onto different planes as fractional-linear mappings of one plane onto another in terms of coordinate systems.  For instance in a supplementary article in Poncelet's book \cite{poncelet1862}, written by one of his collaborators (M. Moutard, presumably a student), one finds the diagram on p. 512 (see Figure \ref{fig:poncelet1862-p512}) which shows a typical projection onto two different planes with coordinate systems \((x,y)\) on the one plane and \((x',y')\) on the second plane.  After four pages of calculations the author finds the coordinate transformations on p. 516 as given in Figure \ref{fig:poncelet1862-p516}.

This use of coordinate systems was used in this context as a way of making sure the computations done via projections agreed with the results being obtained by M\"{o}bius, Pl\"{u}cker and others via coordinate systems.  In fact, the title of the article by Moutard (see p. 509 of \cite{poncelet1862} is: {\it Rapprochement divers entre le principales m\'{e}thodes de la g\'{e}om\'{e}trie pure et celles de l'analyse alg\'{e}brique}.  Here ``{\it rapprochement}"  means reconciliation, and note the use of ``{\it g\'{e}om\'{e}trie pure}" for what later became known as synthetic projective geometry.

\begin{figure}
\vspace{6pt}
\centerline{
	\includegraphics[width=12cm]{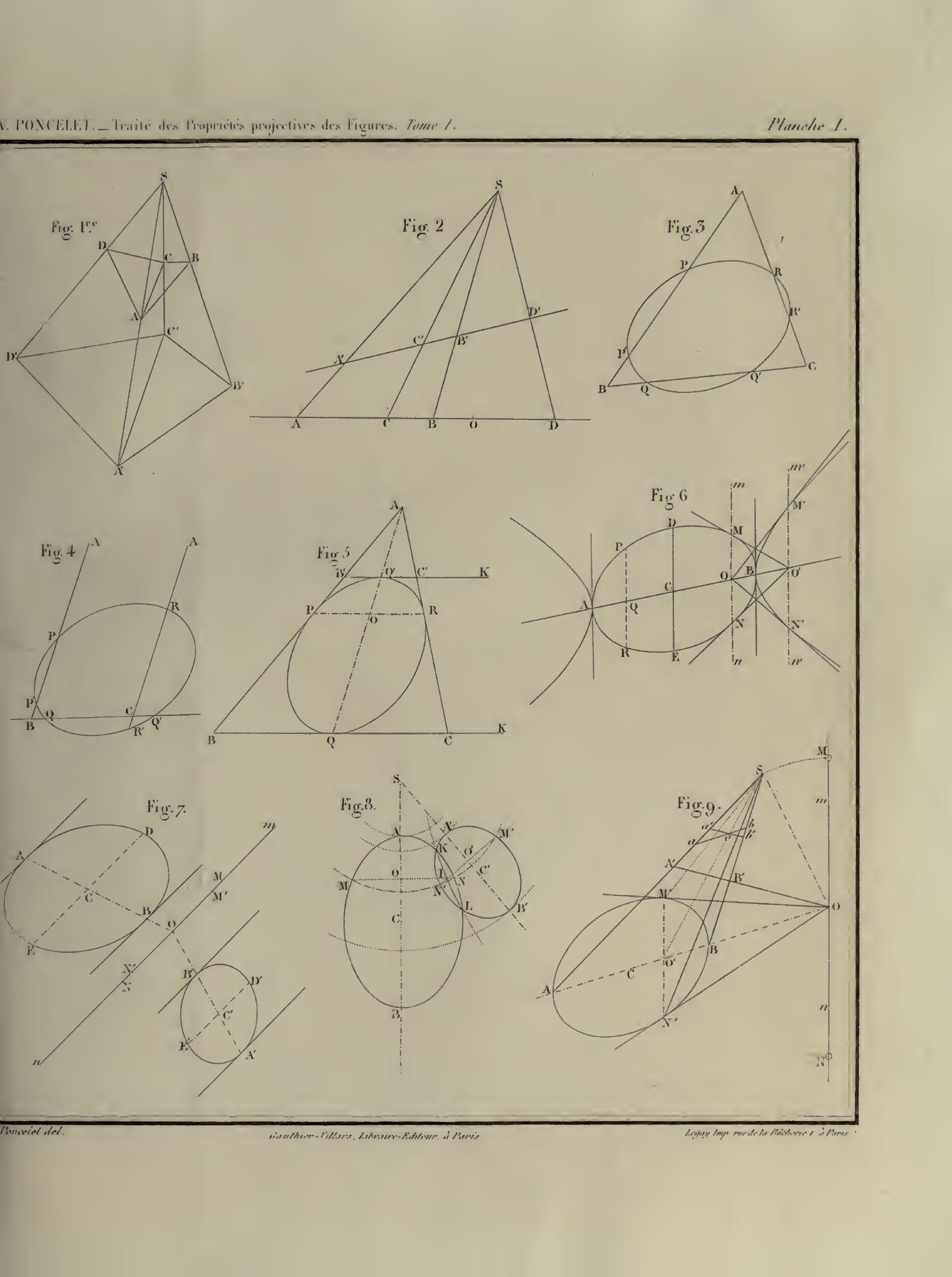}}
 \caption{Figures on Plate 1 in Poncelet's {\it Trait\'{e} des Propri\'{e}t\'{e}s Projectives des Figures, Tome 1.} \cite{poncelet1865}.}
 \label{fig:poncelet1865-fig1}
\end{figure}

The equivalence of geometric objects as being different projections of the same object is one of the essential points of Poncelet's geometry. In his introductory remarks in his books, he formulates the basic principle that guides his investigations: given a particular geometric configuration and an associated problem in this context (often in a planar context), find the simplest projective representation of this figure, by means of which one can resolve the problem, and thus it becomes resolved for all other configurations, including the one in the initial context. We want to illustrate this principle by looking at two simple examples of projectively equivalent figures from Poncelet's book \cite{poncelet1865} (Fig. 1 and Fig. 2 on plate 1, reproduced here in Figure \ref{fig:poncelet1865-fig1}).

In the first figure (Fig. 1) we see that the quadrilateral (assumed in a plane) \(ABCD\) is projectively equivalent to the quadrilateral \(A'B'C'D'\) (again in a plane). In Fig. 2 the segmented line \(ABCD\) is projectively equivalent to the segmented line \(A'B'C'D'\).  Moreover, Poncelet (and many others in the time period) proved that the {\it cross-ratios} of these two sets of collinear points satisfy
\be
\label{eqn:cross-ratio}
\frac{AB}{CB}\dot\frac{CD}{AD}= \frac{A'B'}{C'B'}\dot\frac{C'D'}{A'D'}.
\ee
That is, this numerical quantity doesn't change for these variable projective views of this segmented interval.  It turns out that the cross-ratio is one of the most important numerical invariants of projective geometry. It was first described as lengths of intervals as in (\ref{eqn:cross-ratio}) (and was used by the Greeks in their mathematics work; see the book by Milne \cite{milne1911} for a history of its use in geometry .  Initially it was simply the ratio of lengths of line segments and later was extended to lengths with a sign (when a line was assigned a direction and the order of the points determined the sign, as was introduced by Carnot \cite{carnot1803}). 
This cross-ratio in this context was called in French {\it rapport anharmonique} and in German {\it Doppelverh\"{a}ltnis}, and in modern English we say {\it cross-ratio}, which terminology was introduced by Clifford in \cite{clifford1878} in his study of mechanics. Most contemporary students of mathematics first come across the term cross-ratio in a first course on complex analysis, where for four points \(z_1, z_2, z_3, z_4\) on the extended complex plane \(\overline{\BC} (= \BC \cup \infty)\), the cross-ratio is defined by%
\footnote{Note that there are 24 permutations of these four symbols, and each is called a cross-ratio and satisfies the properties outlined here.  For each such permutation, there are three others with the same values, and hence there are 6 cross-ratios with distinct values. For more details see for instance the very informative book by Milne \cite{milne1911}.}
\[
\frac{z_1-z_2}{z_1-z_4}\dot\frac{z_3-z_4}{z_3-z_2}.
\]
This cross-ratio of four points in the extended complex plane is invariant under fractional-linear transformations (M\"{o}bius transformations) of the form
\[
w=\frac{az+b}{cz+d},
\]
as one learns in a first course in complex analysis (see, e.g., Ahlfors \cite{ahlfors1953}). These fractional-linear transformations are simply reformulations of the projective-linear transformations of the projective space \(\BP_1(\C) \cong \overline{\BC}\).

Poncelet used the word {\it homologous} to describe figures that were projectively equivalent in the sense that we are using here.%
\footnote{Note the similiarity to the modern use of the word``homologous'' in topology, which was introduced by Poincar\'{e} \cite{poincare1895}. In both cases these mathematicians wanted to describe two objects as being similar in their topological or graphical aspects, but not necessarily in their metric aspects ({\it \`{a} la} ``congruence'' or later ``isometry'').}
A second fundamental principle of Poncelet was that of the {\it loi de continuit\'{e}} (law of continuity).  
On p. xiii of his introduction in his basic monograph \cite{poncelet1865} he formulates this principle as an ``axiom" used by illustrious mathematicians of the past. In fact he formulates it as a question, as we see in this excerpt from this page:
\begin{quote}
Consid\'{e}rons une figure quelconque, dans une position g\'{e}n\'{e}rale et en
quelque sorte ind\'{e}termin\'{e}e, parmi toutes celles qu'elle peut prendre sans
violer les lois, les conditions, la liaison qui subsistent entre les diverses parties
du syst\`{e}me; supposons que, d'apr\`{e}s ces donn\'{e}es, on ait trouv\'{e} une ou
plusieurs relations ou propri\'{e}t\'{e}s, soit m\'{e}triques, soit descriptives, appartenant
à la figure, en s'appuyant sur le raisonnement explicite ordinaire, c'est-\`{a} dire
par cette marche que, dans certains cas, on regarde comme seule rigoureuse.
N'est-il pas \'{e}vident que si, en conservant ces m\^{e}mes donn\'{e}es, on
vient \`{a}à faire varier la figure primitive par degr\'{e}s insensibles, ou qu'on imprime
\`{a} certaines parties de cette figure un mouvement continu d'ailleurs
quelconque, n'est-il pas \'{e}vident que les propri\'{e}t\'{e}s et les relations, trouv\'{e}es
pour le premier syst\`{e}me, demeureront applicables aux \'{e}tats successifs de ce
syst\`{e}me, pourvu toutefois qu'on ait \'{e}gard aux modifications particuli\`{e}res
qui auront pu y survenir, comme lorsque certaines grandeurs se seront évanouies,
auront chang\'{e} de sens ou de signe, etc., modifications qu'il sera
toujours ais\'{e} de reconna\^{\i}tre \`{a} priori, et par des r\`{e}gles s\^{u}res?%
\footnote{Consider an arbitrary figure in a general positon and in some sense indeterminant, among all those that one can take
without violating the laws, the conditions, the relations that takes place among the various parts of the system; suppose that, according to this given data, one can find one or more relations or properties, either metric or descriptive, pertaining to the figure, using ordinary expilcit reasoning, that is to say by those steps that, in certain cases one regards as completely rigorous. Isn't is evident that if, in preserving the same given data, one can make the figure vary by imperceptible degrees, or where  one imposes to certain parts of that figure a movement that is continuous and moreover  arbitrary, isn't it evident that the properties and the relations, found for the first system, remain applicable to the successive states of the system, provided that one takes regard of particular modifications that could arise when certain quantities vanish, having a change of directions or sign, etc., modifications that will always be easy to recognize {\it \`{a} priori}, and according to well-determined rules?}
\end{quote}
This quote is indeed his definition of his law of continutiy in this book, and from our point of view, we would say that this is a very {\it vague} definition (to be kind to M. Poncelet), but basically these 19th-century authors used this principle in very many particular cases to derive new results which were later established by more rigorous means.  There were indeed critics of this at the time, e.g. Cauchy was such a critic, and he, during the same period of time, was primarily responsible for formulating many of our current theory of continuity for functions and mappings in general.

\begin{figure}
\vspace{6pt}
\centerline{
	\includegraphics[width=12cm]{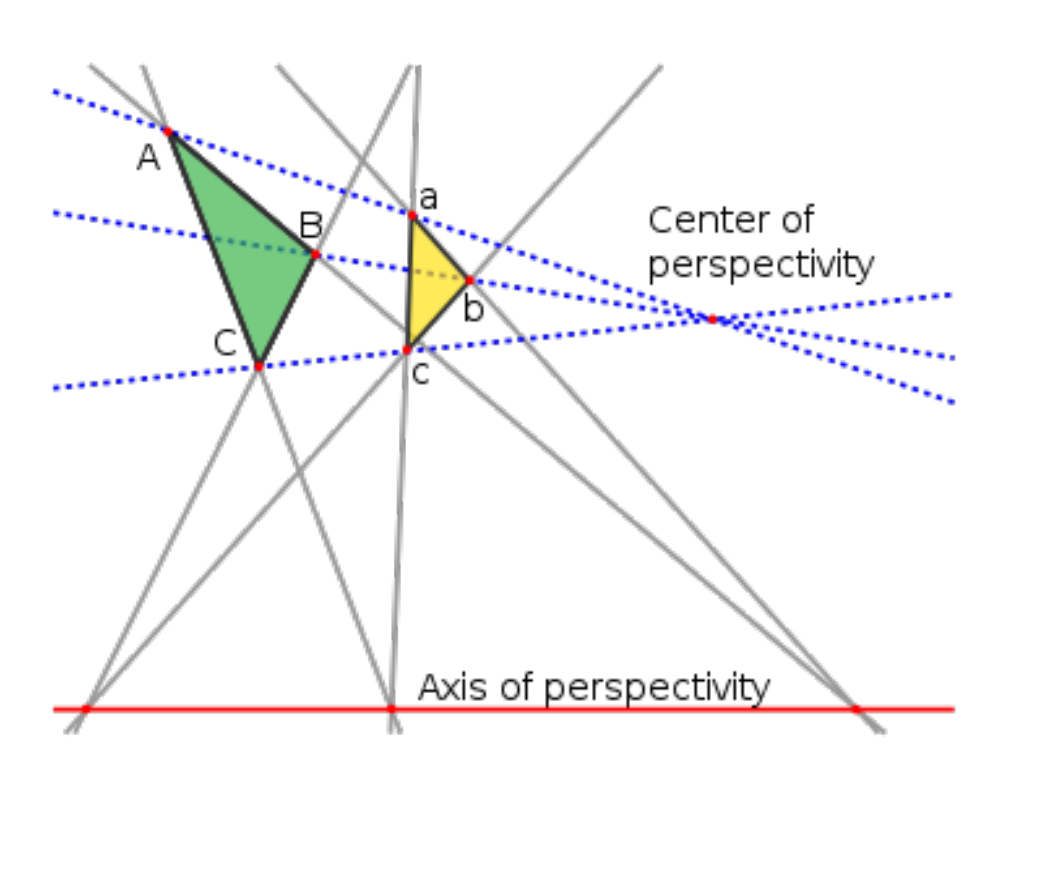}}
 \caption{Illustration of Desargues Theorem.}
 \label{fig:desargues-theorem}
\end{figure}

\begin{figure}
\vspace{6pt}
\centerline{
	\includegraphics[width=12cm]{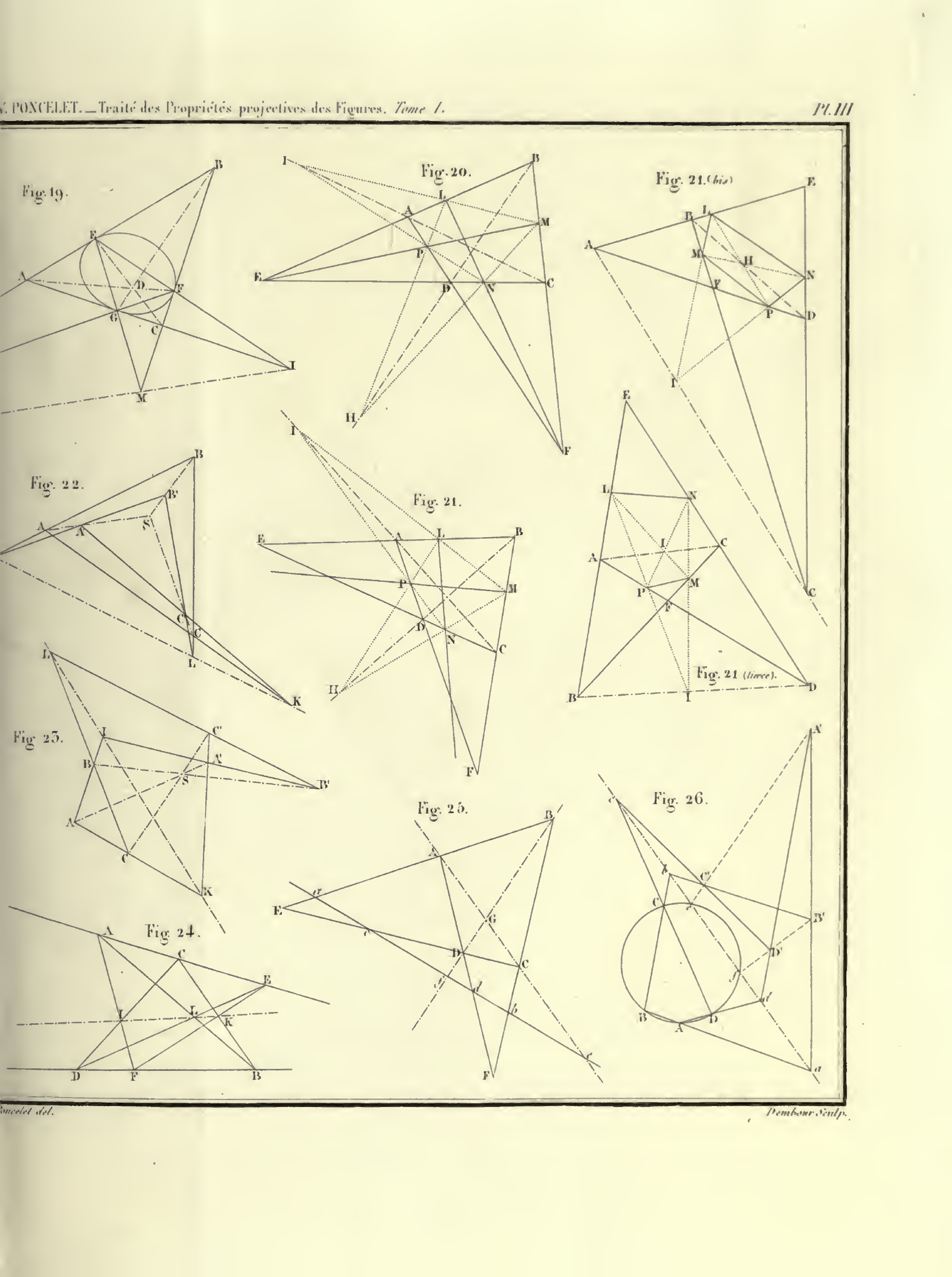}}
 \caption{ Figures on Plate III in Poncelet's {\it Trait\'{e} des Propri\'{e}t\'{e}s Projectives des Figures, Tome 1.} \cite{poncelet1865} .}
 \label{fig:fig22-poncelet1865}
\end{figure}

We want to give an important historical example which illustrates both projective equivalence and the law of continuity, and this is the well-known Desargues Theorem.  In Figure \ref{fig:desargues-theorem}
we see an illustration of the theorem in the plane. Namely, if the two triangles are in perspective with respect to a perspective point (in this case the ``center of perspectivety" in Figure \ref{fig:desargues-theorem}, then the intersections of the homologous line segments lie on a straight line (called the axis of perspective in the figure).  This is the Desargues theorem, and the converse is also true. Let us use Fig. 22 of Poncelet's book (\cite{poncelet1865}, reproduced in Figure \ref{fig:fig22-poncelet1865})
to illustrate the proof.  In this figure we see that the two triangles \(ABC\) and \(A'B'C'\) are in perspective with respect to the point \(S\) and we assume that they are located in two non-parallel planes \(E\) and \(E'\) in \(\BR^3\).  Since both of the lines represented by \(AB\) and \(A'B'\) lie in the plane represented by trangle \(ASB\), then they must necessarily intersect (and for simplicity we assume the intersection is not at infinity). This intersection is denoted by \(M\) (and it is not pictured in Fig. 22 in Figure \ref{fig:fig22-poncelet1865}, since the scan of this page cut off some of the left-hand side of the page, but the intersection is clear on this page).  This intersection point \(M\) lies on the intersection of the two planes \(E\cap E'\), since \(AB\subset E\) and \(A'B' \subset\ E'\).  The same is true for the intersections  \(L\) of the lines represented by \(BC\) and \(B'C'\) and for the intersection \(K\)of the lines represented by \(AC\) and \(A'C'\).  Thus \(K\), \(L\), and \(M\) all lie on the intersection \(E\cap E'\), which is a straight line. This is simply Desargues Theorem in this three-dimensional setting.  Now suppose we have two triangles in perspective in a plane, and we envision them as being continuous limits of triangles in perspective in three-dimensional space, not in the same plane, as above, then the limit of the axes of perspectivity for the three-dimensional case will yield the desired axis of perspectivity in the two-dimensional case.

This fundamental result of projective geometry is due to Desargues in the 17th century. Un fortunately, his work was lost.  However fragments, references to his work, and reworking of some of his work by others was discovered in the course of the beginning of the 19th century, after this theorem, in particular, had been rediscovered by the projective geometry school. One work by Desargues is available today \cite{desargues1642}, which is a draft of three small articles, the first dealing with geometry.  The basic ideas of Desargues are contained in a series of books published from 1643 to 1648 by Bosse, who was a student of Desargues learning about architecture and engraving. We cite the last of this series here \cite{bosse1647} as it is the most mathematical and contains at the end of the book a {\it proposition g\'{e}om\'{e}trique} which is the Desargues Theorem as discussed above. Poncelet gives great credit to Desasrgues and Pascal, whose work we discuss below, for having the initial ideas in projective geometry (see, in particular, p. xxv in the Introduction in \cite{poncelet1865}).

\begin{figure}
\vspace{6pt}
\centerline{
	\includegraphics[width=12cm]{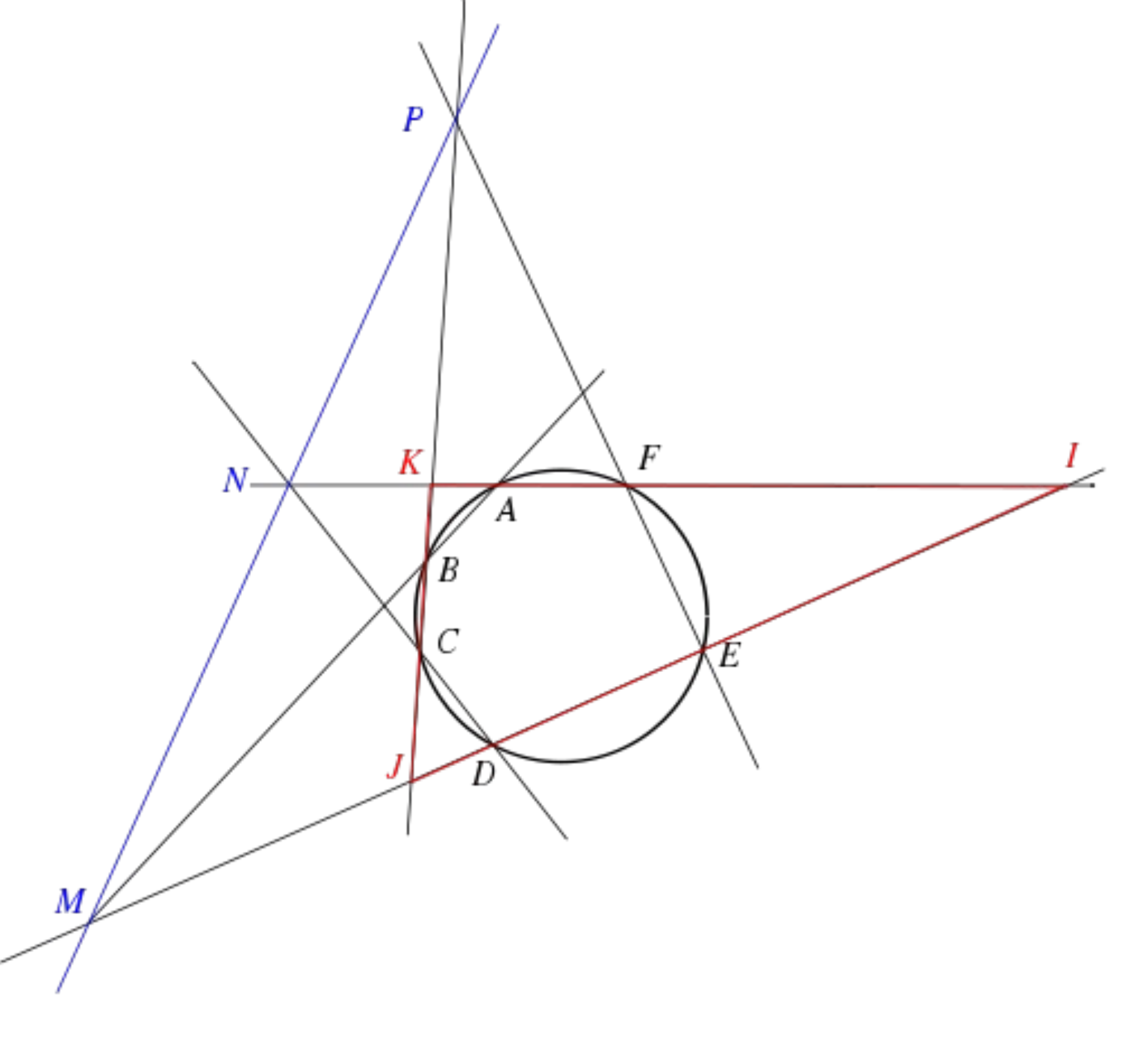}}
 \caption{Ilustration of Pascal's Theorem.}
 \label{fig:pascal-theorem}
\end{figure}

\begin{figure}
\vspace{6pt}
\centerline{
	\includegraphics[width=12cm]{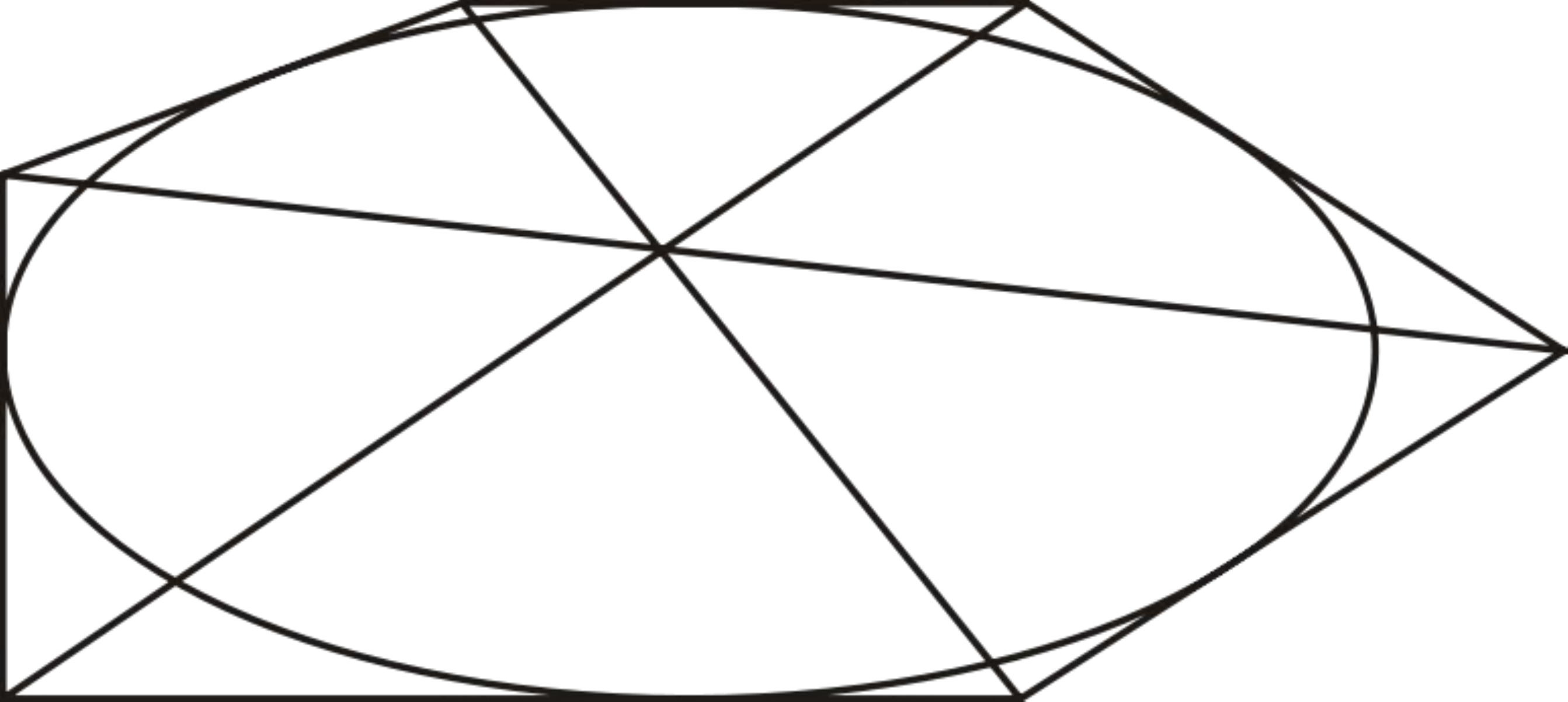}}
 \caption{Illustration of Brianchon's Theorem.}
 \label{fig:brianchon-theorem}
\end{figure}

The final principle of projective geometry we want to mention is that of {\it duality}.  This was first introduced formally in the work of Carnot in 1803 \cite{carnot1803}.  In the simplest case in the plane, by using a nondegenerate quadratic form, one can associate in a one-to-one manner points to lines and conversely,  two points determine a line, and two lines intersect in a point (adding points at infinity here). Propositions utilizing these concepts have dual formulations.  A very classical example of this is illustrated by Pascal's theorem and its dual formulation due to   Brianchon.
Pascal's Theorem asserts that if we consider a hexagram inscribed in a conic section on a plane then the intersections of the opposite sides of the hexagon line on a straight line (see Figure \ref{fig:pascal-theorem} for an illustration of the theorem).
The first published reference to Pascal's theorem is in his collected works, Vol. 5, \cite{pascal1819}, published in 1819.  This volume contains a letter of Leibniz written in 1676 (pp. 429-431 in \cite{pascal1819}) discussing Pascal's papers which he had access to, and including the statement and a proof of Pascal's theorem.  Pascal had written some years earlier a set of essays on conic sections that Leibniz had access to and which he wrote about in this same letter.%
\footnote{The editor (unnamed) of the volume \cite{pascal1819} said that he searched for these papers referred to by Leibniz (written in the mid 17th century) and was not able to find them for this volume of the collected works in 1819.}
 The dual of Pascal's theorem is the theorem of Brianchon (Charles-Julien Brianchon (1783-1864)) \cite{brianchon1806}, which asserts that a circumbscribing hexagon to a conic section has diagonals that intersect at a point (see Figure \ref{fig:brianchon-theorem} for an illustration of this theorem).  Brianchon's monograph \cite{brianchon1817} published a decade later discusses this and related theorems and gives a succinct history of the subject, including references to the work of Pascal (the reference here is simply the title of the paper discussed by Leibniz) and Desargues (there was no mention of Pascal in his paper \cite{brianchon1806}).
Note that the tangent lines in Figure \ref{fig:brianchon-theorem} correspond to the points on the conic section in Pascal's setting, and the intersections of the diagonals, a {\it point}, corresponds to the {\it line} on which the intersections of opposite sides in Pascal's setting lie.

The  final point I'd like to discuss with respect to projective geometry is the use of coordinate systems. In the first decades of the 19th century two schools of projective geometry developed in parallel.  The first school, championed by Poncelet and what we might call the {\it French school} developed what became known as {\it synthetic projetive geometry}, and it was a natural extension of classicasl Euclidean geometry in it methodology, notation, etc. (with the added emphasis on the new ideas in projective geometry, as outlined in the paragraphs above).  The second school was interested in using coordinate systems and analytic notation to prove the essential projective-geometric theorems.  This was called {\it analytic projective geometry} (in the spirit of analytic geometry being the alternative to Euclidean geometry in schools).  These two schools complemented each other, competed in a scholarly manner in the academic publications of the time, and learned  from each other.  The primary authors in the analytic schools happened to be German, and this could be called the {\it German school}. We will discuss several prominent German contributions in the following paragraphs, and we have mentioned the primary participants in the French school earlier.

M\"{o}bius  is the first of the German school we want to mention. In his book {\it Der barycentrische Calcul} \cite{moebius1827} M\"{o}bius was very interested in showing that coordinate systems could be useful in proving the primary theorems of projective geometry (as well as some new results of his own).  He saw the difficulty of using only the standard \((x,y,z)\) coordinate system, which could often be too cumersome or not as elegant as possible (one of the complaints of the French school). M\"{o}bius developed a new coordinate-type system which in its details is very much like the representation of vectors in a vector space using linear combinations of vectors. In \(\BR^3\) he took four points not lying on a plane and took ``linear combinations" of these points (using the notion of center of mass from mechanics).  Using the formalism he developed, he was able to successfully develop the fundamental results in projective geometry.  Today his ideas are still used in barycentric subdivision in the triangulation of topological spaces, and it is a quite original, very readable, and interesting book to peruse. In this book he classifies geometric structures according to {\it congruence} (he uses the word ``{\it Gleichheit}", that is equivalent under Euclidean motions (translations and rotations), {\it affine} that is translations and linear transformations of \(\BR^3\), and {\it collinear} (equivlence under mappings preserving lines, projective equivalence as described earlier in this section).

The second author we turn to, whose work turned out to be very fundamental is Pl\"{u}cker.  He wrote two books {\it Analytisch-geometrische Entwicklungen} Vol. 1 \cite{pluecker1828} in 1828 and Vol. 2 \cite{pluecker1831} in 1831.  In the first volume Pl\"{u}cker developed an abridged analytic notation to make the analytic proofs more transparent  and more concise (abridgement of the classical coordinate system notation in \(\BR^3\)).  In Vol. 2 from 1831, he excitedly told his readers that he had developed yet another and new notation which simpliefied the understnding of projective geometry, and this was his introduction of {\it homogeneous coordinates}, which he had announced in a research paper in Crelle's Journal in 1830 \cite{pluecker1830}. Later, after a successful career in the spectroscopy of gases, he returned to the study of geoemtry and developed the notion of {\it line geometry}, the study of the lines in \(\BR^3\) as basic objects of study, and which could be parametrized as a quadric surface in \(\BP_5\). This work was first developed in a paper published in 1865 \cite{pluecker1865} and developed into a two-part monograph as a new way of looking at geometrical space.  This monograph \cite{pluecker1868} was published after Pl\"{u}cker's death and was edited by A. Clebsch with a great deal of assistance by the young Felix Klein, who was  a student at the University of Bonn studying with Pl\"{u}cker at the time. The fundamental idea was to make lines and their parametrizations the principal object of study (in contrast with having points in space being the primary objects).  One obtained points as intersections of lines, just as in classical geometry one obtained a line passing through two points. This led naturally (among many other things) to the definition of  two-dimensional real projective space as the space of all lines passing through a given point in \(\BR^3\). 

The most original contribution to the analytic side of projective geometry (indeed to geometry itself!) in this time period was that of Grassmann  \cite{grassmann1844}, who, in a singular work, formulated what became know today as {\it exterior algebra}.  His work also laid the foundation for the theory of vector spaces, one of the building blocks for exterior algebra. His work led to the development of Grassmannian manifolds, which regarded all of the planes of a fixed dimension in a Euclidean space as a geometric space, which was a generalization of the line geometry of Pl\"{u}cker.
His lengthy philosophical introduction to his work is both a challenge and at the same time a pleasure to read.  His work was not that well recognized in his lifetime (as so often happened in the history of mathematics), but has turned out in the hands of Elie Cartan, for instance, to be a powerful tool for studying geometric problems in the 20th century.

Finally we note that Felix Klein in his writings put projective geometry on the firm footing we see today.  In particular his {\it Erlangen Program}\footnote{This was written as his research program when he took up a new professorship in Erlangen in 1872.} \cite{klein1872} from 1872 and his lectures on non-Euclidean geometry from the 1890s, published posthumously in 1928 \cite{klein1928}, are marvels of exposition and give us, for instance the use of the term {\it projective space} as a concept. 

\section{Riemann's Higher-Dimensional Geometry}
In mathematics we sometimes see striking examples of brilliant contributions or completely new ideas that change the ways mathematics develops in a significant fashion.   A prime example of this is the work of Descartes \cite{descartes}, which completely changed how mathematicians looked at geometric problems.  But it is rare that a single mathematician makes as many singular advances in his lifetime as did Riemann in the middle of the 19th century. In this section we will discuss in some detail his fundamental creation of the theory of higher-dimensional manifolds and the additional creation of what is now called Riemannian or simply differential geometry.  However, it is worth noting that he only published nine papers in his short lifetime (he lived to be only 40 years old), and several other important works, including those that concern us in this section, were published posthumously from the writings he left behind.  His collected works (including in particular these posthumously published papers) were edited and published in 1876 and are still in print today \cite{riemann1876}.
\begin{figure}
\vspace{6pt}
\centerline{
	\includegraphics[width=12cm]{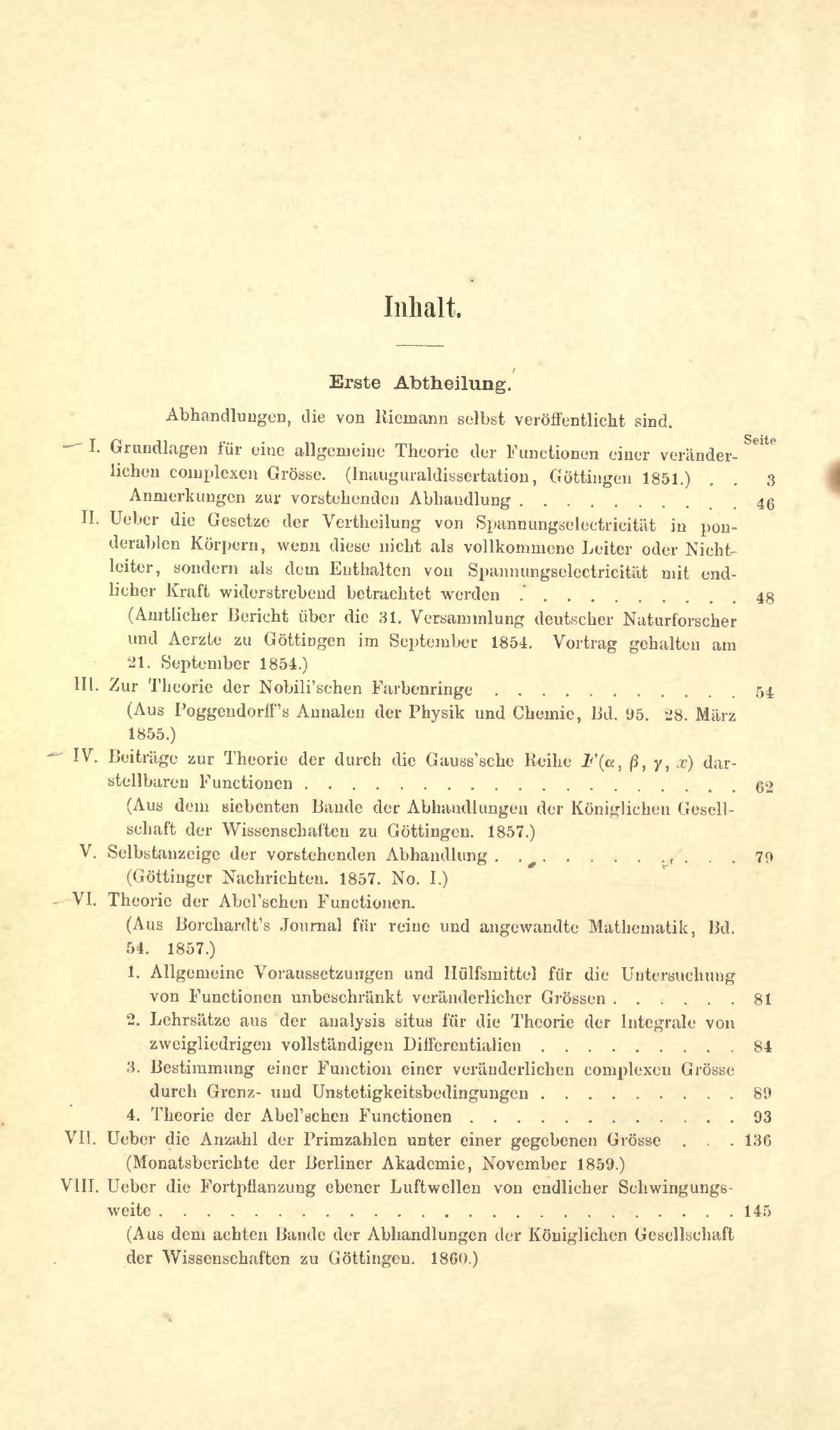}}
 \caption{Table of Contents p. VI, Riemann's Collected Works \cite{riemann1876}.}
 \label{fig:riemann-TOC-1}
\end{figure}
\begin{figure}
\vspace{6pt}
\centerline{
	\includegraphics[width=12cm]{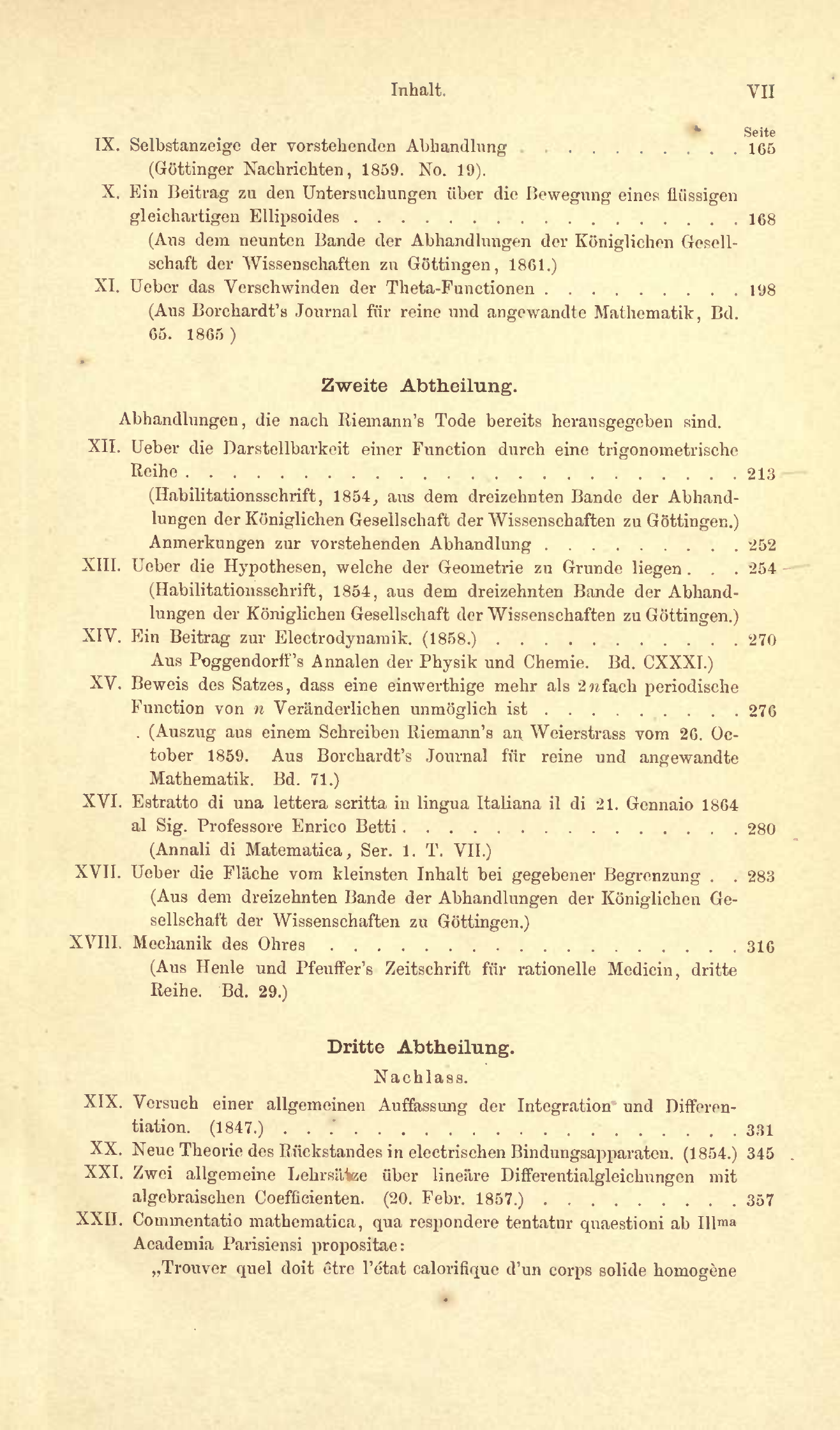}}
 \caption{Table of Contents p. VII, Riemann's Collected Works \cite{riemann1876}.}
 \label{fig:riemann-TOC-2}
\end{figure}

Looking through the titles one is struck by the wide diversity as well as the originality.  Let us give a few examples here.  In Paper I (his dissertation) he formulated and proved the Riemann mapping theorem and dramatically moved the theory of functions of one complex variables in new directions.  In Paper VI, in order to study Abelian functions, he formulated what became known as Riemann surfaces and this led to the general theory of complex manifolds in the 20th century.  In Paper VII he proved the Prime Number Theorem and formulated the Riemann Hypothesis, which is surely the outstanding mathematical problem in the world today. In Paper XII he formulated the first rigorous definition of a definite integral (the Riemann integral) and applied it to trigonometric series, setting the stage for Lebesgue and others in the early twentieth century to develop many consequences of the powerful theory of Fourier analysis.  In Papers XIII and XXII he formulated the theory of higher-dimensional manifolds, including the important concepts of Riemannian metric, normal coordinates and the Riemann curvature tensor, which we will visit very soon in the paragraphs below.  Paper XVI contains correspondence with Enric Betti leading to the first higher-dimensional topological invariants beyond those already known for two-dimensional manifolds.

This will suffice.  The reader can glance at the other titles to see their further diversity.  His contributions to the theory of partial differential equations and various problems in mathematical physics were also quite significant.

His paper \cite{riemann1854} (Paper XIII above) is a posthumously published version of a public lecture Riemann gave as his { \it Habilitationsvortrag} in 1854. This was part of the process for obtaining his {\it Habilitation}, a German advanced degree beyond the doctorate necessary to qualify for a professorship in Germany at the time (such requirements still are in place at most German universities today as well as in other European countries, e.g., France and Russia; it is similar to the research requirements in the US to be qualified for tenure). This paper, being a public lecture, has very few formulas, is at times quite philosophical and is amazing in its depth of vision and clarity.  On the other hand, it is quite a difficult paper to understand in detail, as we shall see.

Before this paper was written, manifolds were all one- or two-dimensional curves and surfaces in \(\BR^3\), including their extension to points at infinity, as discussed in Section \ref{sec:projective geometry}.  In fact some mathematicians who had to study systems parametrized by more than three variables declined to call the parametrization space a manifold or give such a parametrization a geometric significance.  In addition, these one- and two-dimensional manifolds always had a differential geometric structure which was induced by the ambient Euclidean space (this was true for Gauss, as well).

In Riemann's paper \cite{riemann1854} he discusses the distinction between discrete and continuous manifolds, where one can make comparisons of quantities by either counting or by measurement, and gives a hint, on p. 256, of the concepts of set theory, which was only developed later in a single-handed effort by Cantor). Riemann begins his discussion of manifolds by moving a one-dimensional manifold, which he intuitively describes, in a transverse direction (moving in some type of undescribed ambient "space"), and inductively, generating an \(n\)-dimensional manifold by moving an \(n-1\)-manifold transversally in the same manner.  Conversely, he discusses having a nonconstant function on an \(n\)-dimensional manifold, and the set of points where the function is constant is (generically) a lower-dimensional manifold; and by varying the constant, one obtains a one-dimensional family of \(n-1\)-manifolds (similar to his construction above).%
\footnote{He alludes to some manifolds that cannot be described by a finite number of parameters; for instance, the manifold of all functions on a given domain, or all deformations of a spatial figure.  Infinite-dimensional manifolds, such as these, were studied in great detail a century later.}

Riemann formulates local coordinate systems \((x^1, x^2,...,x^n)\) on a manifold of \(n\) dimensions near some given point, taken here to be the origin. He formulates a curve in the manifold as being simply \(n\) functions \((x^1(t), x^2(t),...,x^n(t))\) of a single variable \(t\). The concepts of set theory and topological space were developed only later in the 19th century, and so the global nature of manifolds is not really touched on by Riemann (except in his later work on Riemann surfaces and his correspondence with Betti, mentioned above). It seems clear on reading his paper that he thought of \(n\)-dimensional manifolds as being extended beyond Euclidean space in some manner, but the language for this was not yet available.

At the beginning of this paper Riemann acknowledges the difficulty he faces in formulating his new results.  Here is a quote from the second page of his paper (p. 255):
\begin{quote}
Indem ich nun von diesen Aufgaben zun\"{a}chst die erste, die Entwicklung
des Begriffs mehrfach ausgedehnter Gr\"{o}ssen, zu l\"{o}sen versuche,
glaube ich um so mehr auf eine nachsichtige Beurtheilung Anspruch
machen zu d\"{u}rfen; da ich in dergleichen Arbeiten philosophischer
Natur, wo die Schwierigkeiten mehr in den Begriffen,
als in der Construction liegen, wenig ge\"{u}bt bin und ich ausser einigen ganz kurzen
Andeutungen, welche Herr Geheimer Hofrath Gauss in der zweiten
Abhandlung \"{u}ber die biquadratischen Reste in den G\"{o}ttingenschen
gelehrten Anzeigen und in seiner J{u}bil\"{a}umsschrift dar\"{u}ber gegeben
hat, und einigen philosophischen Untersuchungen Herbart's, durchaus
keine Vorarbeiten benutzen konnte.%
\footnote{In that my first task is to try to develop the concept of a multiply spread out quantity [he uses the word manifold later], I believe even more in  being allowed an indulgent evaluation, as in such works of a philosophical nature, where the difficulties are more in the concepts than in the construction, wherein I have little experience, and except for the paper by Mr. Privy Councilor Gauss in his second commentary on biquadratic residues in the G\"{o}inngen gelehrten Anzeigen [1831]  and in his J{u}bil\"{a}umsschrift and some investigations by Hebart, I have no precedents I could use.}
\end{quote}
The paper of Gauss that he cites here \cite{gauss1831} refers to Gauss's dealing with the philosophical issue of understanding the complex number plane after some thirty years of experience with its development.  We will mention this paper again somewhat later in this paper.  Hebart  was a philosopher whose metaphysical investigations influenced Riemann's thinking.  Riemann was very aware of the speculative nature of his theory, and he used this philosophical point of view, as the technical language he needed (set theory and topological spaces) was not yet available.  This was very similar to Gauss's struggle with the complex plane, as we shall see later.

As mentioned earlier, measurement of the length of curves goes back to the Archimedian study of the length of a circle.  The basic idea there and up to the work of Gauss was to approximate a given curve by a straight line segment and take a limit.  The {\it length} of each straight line segment was determined by the Euclidean ambient space, and the formula, using calculus for the limiting process, became, in the plane for instance,
\[
\int_\G ds = \int^b_a\sqrt{(x'(t))^2 + (y'(t))^2}dt,
\]
where \(ds^2=dx^2+dy^2\) is the line element of arc length in \(\BR^2\).  As we saw in Section \ref{sec:gauss}, Gauss formulated in \cite{gauss1828} on a two-dimensional manifold with coordinates \((p,q)\) the line element
\be
\label{eqn:line-element}
ds^2=Edp^2 + 2Fdpdq + Gdq^2
\ee
where \(E,F,\) and \(G\) are induced from the ambient space.  He didn't consider any examples of such a line element (\ref{eqn:line-element}) that weren't induced from an ambient Euclidean space, but his remarks (see the quote above in Section \ref{sec:gauss}) clearly indicate that this could be a ripe area for study, and this could well include allowing coeffients of the line element  (\ref{eqn:line-element}) being more general than induced from an ambient space.

Since Riemann formulated an abstract \(n\)-dimensional manifold (with a local coordinate system) with no ambient space, and since he wanted to be able to measure the length of a curve on his manifold, he formulated, or rather postulated, an independent measuring system which mimics Gauss'ss formula (\ref{eqn:line-element}).  Namely, he prescribes for a given coordinate system a metric (line element) of the form 
\be
\label{eqn:riemann-metric}
ds^2= \sum^n_{i,j=1} g_{ij}(x)dx^i dx^j,
\ee
where  \(g_{ij}(x) \) is, for each \(x\), a positive definite matrix, and he postulates by the usual change of variables formula
\[
ds^2= \sum^n_{i,j=1}\tilde{g}_{ij}(\tilde{x}) d \tilde{x}^i d \tilde{x}^j,
\]
where \(\tilde{g}_{ij}(\tilde{x})\) is the transformed positive definite matrix in the new coordinate system \((\tilde{x}_1,...,\tilde{x}^n
)\). 

Using the line element (\ref{eqn:line-element}), the length of a curve is defined by 
\[
l(\G):= \int^b_a\sqrt{\sum^n_{i,j=1} g_{ij}(x(t))\frac{dx^i}{dt}(t)\frac{dx^j}{dt}(t)}dt.
\]
The line element (\ref{eqn:riemann-metric}) is what is called a {\it Riemannian metric} today, and the 2-form \(ds^2\) is considered as a positive definite bilinear form giving an inner product on the tangent plane \(T_p (M)\) for \(p\) a point on the manifold \(M\). This has become the basis for almost all of modern differential geometry (with the extension to Lorentzian type spaces where \(g_{ij}(x)\) is not positive definite {\it \`{a} la} Minkoswki space, or symplectic geometry, etc.). Riemann merely says
on page 260 of his paper (no notation here at all),
\begin{quote}
ich beschr\"{a}nke mich daher auf die
Mannigfaltigkeiten, wo das Linienelement durch die Quadratwurzel aus
einem Differentialausdruck zweiten Grades ausgedr\"{u}ckt wird.%
\footnote{I restrict myself therefore to manifolds where the line element is expressed by the square root of a differential expression of second degree.}
\end{quote}
Earlier he had remarked that a line element should be homogenous of degree 1 and one could also consider the fourth root of a differential exression of fourth degree, for instance.  Hence his restriction in the quote above.

\begin{figure}
\vspace{6pt}
\centerline{
	\includegraphics[width=12cm]{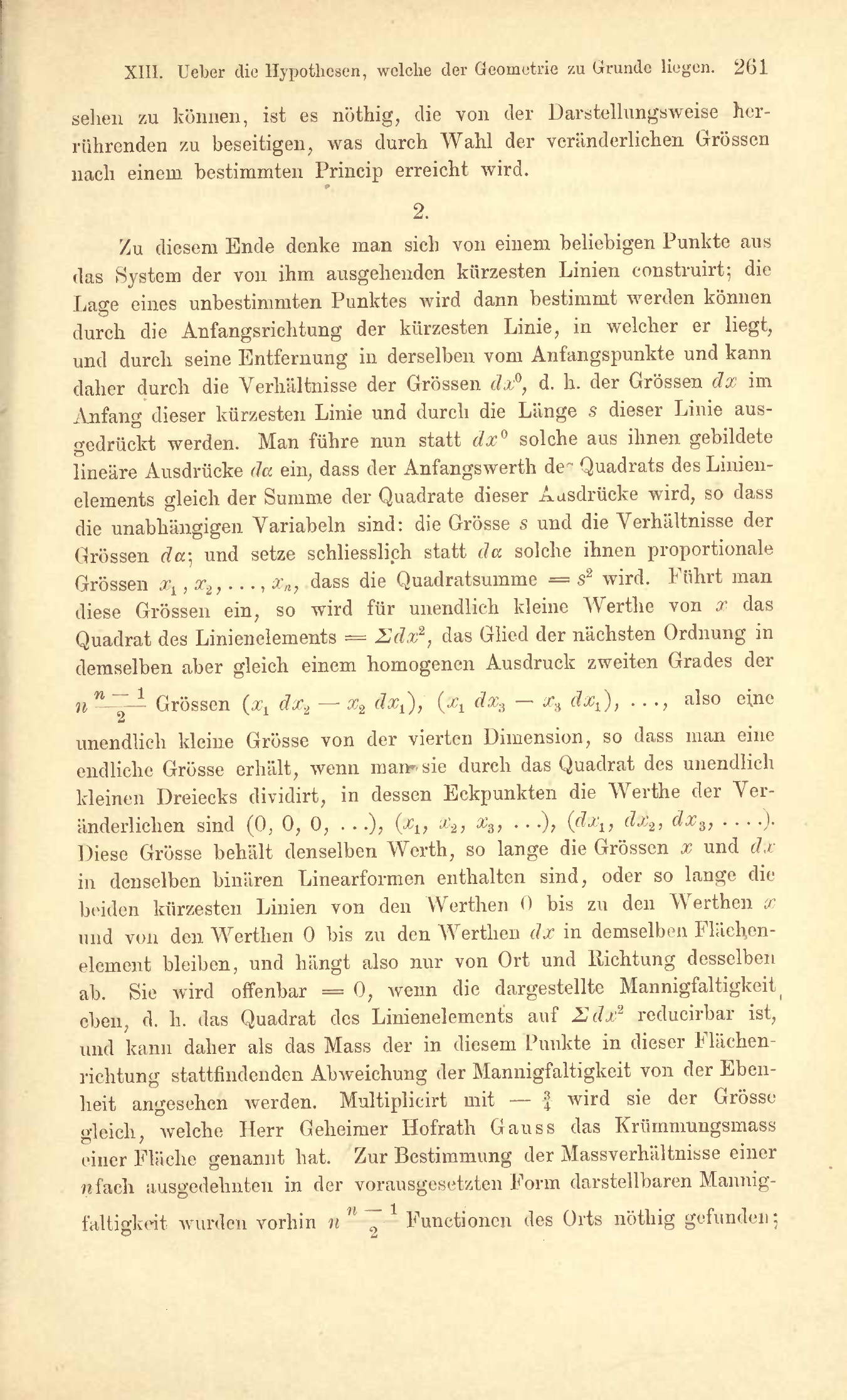}}
 \caption{Page 261 of Riemann's  foundational paper on differential geometry \cite{riemann1854}.}
 \label{fig:riemann1854p261}
\end{figure}

The next step in Riemann's paper is his formulation of curvature.  This occurs on a single page (p. 261 of \cite{riemann1854}, which we reproduce here in Figure \ref{fig:riemann1854p261}).
It is extremely dense and not at all easy to understand.  In the published collected works of Riemann one finds an addendum to Riemann's paper which analyzes this one page in seven pages of computations written by Dedekind.  This is an unpublished manuscript published only in these collected works of Riemann, pp. 384-391). In Volume 2 of Spivak's three volume comprehensive introduction to and history of differential geometry \cite{spivak1999-2}, we find a detailed analysis of Riemann's paper (as well as Gauss's papers that we discussed earlier and later important papers of the 19th century in differential geometry, including translations into English of the most important papers).

We want to summarize what Riemann says on p. 261 (again, see Figure \ref{fig:riemann1854p261}). He starts by introducing near a given point \(p\) on his manifold \(M\) {\it geodesic normal coordinates}, that is, coordinates which are geodesics emanating from the given point and whose tangent vectors at \(p\) are an orthonormal basis for \(T_p\)
(this orthogonality and the geodesics use, of course, the given Riemannian metric). In this coordinate system \((x^1,...,x^n)\), the metric \(ds^2\) has a Taylor expansion through second-order terms of the form
\be
\label{eqn:taylor}
ds^2= \sum^n_{i=1}dx^i dx^i + \frac{1}{2}\sum^n_{ijkl} \frac{\partial^2 g_{ij}}{\partial x^k\partial x^l} (0)x^k x^l dx^idx^j .
\ee
The first-order terms in this expansion involve terms of the form \(\frac{\partial g_{ij}}{\partial x^k}(0)\), all of which vanish, which follows from the geodesic coordinates condition.
Letting now
\[
c_{ijkl} := \frac{\partial^2 g_{ij}}{\partial x^k \partial x^l}(0),
\]
we have the natural symmetry conditions
\[
c_{ijkl} = c_{jikl}=c_{jilk},
\]
due to the symmetry of the indices in \(g_{ij}\) and in the commutation of the second-order partial derivatives. Moreover, and this is {\it not} easy to verify, the coefficients {\it also} satisfy
\be
\label{eqn:symmetry}
\begin{array}{c}
c_{ijkl}= c_{klij},\\
c_{lijk} + c_{ljki} + c_{lkij}= 0.
\end{array}
\ee
This is proved in six pages of computation in Spivak's Vol. 2 (pp. 172--178 of \cite{spivak1999-2}, and we quote from the top of p. 174 ``We now proceed to the hardest part of the computation, a hairy computation indeed."
These symmetry conditions use the fact that the coordinates are specifically linked to the metric (our geodesic coordinates).  For instance, on p. 175 Spivak points out that
\[
x^i = \sum^n_{j=1} g_{ij}x^j,
\]
illustrating vividly the relation between the coordinates and the metric.

Let now 
\be
\label{eqn:biquadratic}
Q(x,dx):= \sum_{ijkl}c_{ijkl} x^kx^ldx^idx^j
\ee
be the biquadratic form defined by the second-order terms in (\ref{eqn:taylor}), then Riemann asserts on p. 61 of \cite{riemann1854} that \(Q(x,dx)\) can be expressed in terms of the \(n(\frac{n-1}{2})\) expressions \( \{(x^1dx^2-x^2dx^1),\) \((x^1dx^3-x^3dx^1)\),,...\}. (see Figure \ref{fig:riemann1854p261}), that is, 
\be
\label{eqn:skew-symmetry}
Q(x,dx) = \sum^n_{ijkl} C_{ijkl}(x^idx^j-x^jdx^i)(x^kdx^l-d^ldx^k).
\ee
Spivak proves that the conditions  (\ref{eqn:symmetry}) are necessary and sufficient for \(Q(x,dx)\) to be expressed in the form (\ref{eqn:skew-symmetry}), and he shows moreover, that
\[
C_{ijkl} = \frac{1}{3} c_{ijkl}.
\]
Riemann simply asserts that this is the case, which is, of course, indeed true! 

Let's look at the special case where \(\dim M = 2\).  In this case we see that there is only one coefficient of the nonzero term \((x^1dx^2-x^2dx^1)^2\) which has the form
\[
Q(x,dx) = \frac{1}{3} [ c_{2211} +c_{1122} - c_{2112} -c_{1221}](x^1dx^2-x^2dx^1)^2.
\]
Now using Gauss's notation for the Riemann metric (\ref{eqn:line-element}), that is, \(g_{11}= E\), \(g_{12}=g_{21} = F\), and \(g_{22}= G\), we see that
\bean
c_{2211}	&=	&\frac{1}{2}G_{xx},\\
c_{1122}	&=&	\frac{1}{2}E_{yy},\\
c_{2112} &=& \frac{1}{2} F_{xy}.\\
c_{1221} &=& \frac{1}{2} F_{xy},
\eean
and thus we have
\[
Q(x,dx) = \frac{1}{6}[G_{xx} +E_{yy} -2F_{xy}].  
\]
Looking at Gauss's formula for Gaussian curvature at the point \(0\) (\ref{eqn:egregrium}), we see that, since the first derivatives of the metric vanish at the origin, the curvature at \(x=0\) is
\be
\label{eqn:curvature}
k= -\frac{1}{2}(G_{xx}G_{yy}-2F_{xy}),
\ee
and hence
\be
\label{eqn:q-simple}
Q(x,dx) = -\frac{k}{3}(x^1dx^2-x^2dx^1)^2.
\ee
Thus the coefficient of the single term \((x^1dx^2-x^2dx^1)^2\) in the biquadratic form \(Q(x,dx)\) is, up to a constant, the Gaussian curvature. As Riemann asserts it: divide the expression \(Q(x,dx)\) by the square of the area of the (infinitesimal) triangle  formed by the three points \((0,x,dx)\), and the result of the division is \(-\frac{4}{3}k\).  The factor 4 appears, since the square of the area of the infinitesimal parallelogram%
\footnote{Riemann visualizes the parallelogram formed by the points \((0,x,dx, x+dx)\) in \(\BR^2\) and the area of such a rectangle is simply given by the cross product \(\|x \times dx\|=\|x^1dx^2-x^2dx^2\|\), and the area of the triangle formed \((0,x,dx)\) is \(\frac{1}{2} \|x^1dx^2-x^2dx^1\|\).}
 is \(x^1dx^2-x^2dx^1)^2\)
 and thus the square of the  area of the infinitesmal triangle is \(\frac{1}{4}(x^1dx^2-x^2dx^1)^2\). This yields the relation between Riemann's coefficient in (\ref{eqn:q-simple}) and Gaussian curvature (one can see this coefficient of \(-\frac{3}{4}\) near the bottom of  p. 261 in Figure \ref{fig:riemann1854p261}).
Riemann then considers the biquadratic form \(Q(x,dx)\) in an \(n\)-dimensional manifold \(M\) and its restriction to any two-dimensional submanifold \(N\) passing through the point \(p\), obtaining a curvature (constant multiple of the Gaussian curvature as we saw above) for the submanifold at that point.  This is the sectional curvature of Riemann introduced on this same p. 261.

In the remainder of the paper he discusses questions of flat manifolds, manifolds of positive or negative constant curvature, and numerous other questions.

The coefficients \(\{c_{ijkl}\}\) in (\ref{eqn:biquadratic}) are effectively the components of the Riemannian curvature tensor for this special type of coordinate system (geodesic normal coordinates).  How does one define such a curvature tensor for \(n\)-dimensional manifolds with a Riemannian metric in a general coordinate system (in the spirit of Gauss's curvature formula (\ref{eqn:egregrium}))?  Clearly this will involve the first derivatives of the Riemannian metric as well.   In a  paper written in Latin for a particular mathematical prize in Paris, (Paper No. XXII in Figure \ref{fig:riemann-TOC-2}) Riemann provides the first glimpse of the general Riemann curvature tensor, and this is again translated and elaborated on by Spivak \cite{spivak1999-2}.  The purpose of this paper was to answer a question in the Paris competition dealing with the flow of heat in a homogenous solid body.

Riemann's ideas in these two posthumously published papers were developed and expanded considerably in the following decades in the work of Christoffel, Levi-Cevita, Ricci, Beltrami and many others. This is all discussed very elegantly in Spivak's treatise \cite{spivak1999-2}, and we won't elaborate this any further at this point. The main point of our discussion has been that Riemann created on these few pages the basic idea of an \(n\)-dimensional manifold not considered as a subset of Euclidean space {\it and} of the independent concept of a Riemannian metric and the Riemannian curvature tensor.  What is missing at this point in time is the notion of a topological space on the  basis of which one can formulate the contemporary concept of an abstract manifold or an abtract Riemannian manifold.

\section{Concluding Remarks: On to the 20th Century}
\markboth{Concluding Remarks}{Concluding Remarks}
This concludes our discussion of three of the key developments in geometry in the 19th century.  Other major developments in geometry that led to the study of abstract manifolds in the 20th century, and which we won't consider in any detail in this paper are the creations of complex geometry,  transformations groups, set theory and  topology, among others. Complex numbers were called {\it imaginary numbers} (among other appelations) in the centuries  preceding the 19th century.  This meant, in the eyes of the beholders that these were not concrete real (or realistic) numbers (to use a pun!) but were simply imaginary artifacts that had no real meaning, but were useful in the way they arose as the would-be solutions of algebraic equations. Gauss, among others, used these numbers extensively in his career, for instance in his ground-breaking work on number theory early in his life {\it Disquisitiones Arithmeticae} \cite{gauss1801}, and only many decades later did he make the case for a definitive geometric interpretation of what we call complex numbers today in a brief paper \cite{gauss1831}, which was a commentary on some of his earlier work on number theory and which Riemann cited earlier in this paper as a philosophical work which was a guide for him. Gauss called these numbers (for the first time) {\it complex numbers}, with their real and imaginary parts representing coordinates in a two-dimensional plane --- the {\it complex plane}. Gauss pleads with his readers to consider complex numbers and the complex plane to be considered as a normal part of mathematics, not as something "imaginary" or "unreal".  It's very enlightening to read this very readable short paper.  Later in the century Klein and others started to consider complex solutions of homogenous algebraic equations and thus began the study of complex algebraic manifolds and varieties.  Klein along with Lie fostered the study of transformation groups, and the notion of manifolds as quotients of such groups became an important development as well.  See, for instance, the very informative book by Klein \cite{klein1928} (first published in 1928, but based on lectures of Klein towards the end of the 19th century), which described most of the developments in geometry in the 19th century in a succinct fashion.

The creation of set theory in the latter half of the 19th century was a singular effort of Georg Cantor over several decades of work.  The first article in this direction was published in 1874 \cite{cantor1874} He then published six major papers in the Mathematische Annalen between 1879 and 1884,%
\footnote{The author  recalls a discussion many years ago with Hans Grauert, at the time Managing Editor of the {\it Mathematische Annalen}, who mentioned that he had access to all of the editorial files of his predecessors, and that he had learned from this corespondence how very controversial it had been at that time to publish these papers.}
which established the basic tenets of set theory, and laid the foundation for work in multiple directions for the next century and beyond, in particular in logic and foundations of mathematics, as well as point set topology, to mention a particular part of geometry related to our thesis in this paper (see his collected works \cite{cantor1932} for these papers, as well as numerous reference works on basic set theory, which every undergraduate learns today).

After  Cantor's goundbreaking discovery (which was very controversial for quite some time) the notion of topological space evolved to its contemporary state as a set satisfying certain axioms concerning either neighborhoods or open sets.  This became the development of point set topology, which has its own very interesting history, and we mention only the fundamental contributions of Frechet in 1906 \cite{frechet1906} (metric spaces, Frechtet spaces, in particular in the infinite dimensional setting) and Hausdorff in 1914 \cite{hausdorff1914} (which, among other things, defined specifically Hausdorff spaces, one of many specialized types of topological spaces used often in geometry today). In 1895 Poincar\'{e} launched the theory of topological manifolds with a cornerstone paper called {\it Analysis Situs} \cite{poincare1895}, which was followed up by five supplements to this work over the next 17 years (see Vol. 6 of his collected works \cite{poincare1953} as well as a very nice translation of all of these topology papers of Poincar\'{e} by John Stillwell \cite{poincare2010}). This paper and its supplements became the foundation of what is now called {\it algebraic topology} (following the pioneering work of Riemann \cite{riemann1857} on Riemann surfaces and Betti on homological invariants of higher dimensional manifolds \cite{betti1871}).

The final step in our journey  is the book by Hermann Weyl in 1913 \cite{weyl1913} entitled {\it Die Idee einer Riemannsche Fl\"{a}che},%
\footnote{The Concept of a Riemann Surface.}
in which he gives the first very specific definition of an abstract manifold with more structure than a topological manifold (which Poincar\'{e} and others had already investigated quite thoroughly).  His motivation was to give a better understanding of Riemann surfaces which transcended the picture of a multisheeted covering of the complex plane with branch points of various kinds. Specifically he considered a topological manifold of two dimensions (locally homeomorphic to an open set in \(\BR^2\) ) which had a finite or countable triangulation and which had the additional property that there were coordinate charts mapping to an open ball in the complex plane \(\BC\) whose transition functions on overlapping coordinate charts were holomorphic. He showed how all of the previous work on Riemann surfaces fit into this new picture, and he proved a fundamental existence theorem of global holomorphic or meromorphic functions on such surfaces, utilizing the Dirichlet principle, that had been first used by Riemann in his conformal mapping theorem from his dissertation \cite{riemann1851}.  This definition of a topological manifold with such additional structure became the model for all the various kinds of manifolds studied in the following century up to the current time. This included, for instance, differentiable (\(C^\infty\) ) manifolds, complex manifolds of arbitrary dimension, Riemannian manifolds, symplectic manifolds, real-analytic manifolds, among many others. Weyl knew that this was new territory, and like Gauss with intrinsic differential geometry, and Riemann with \(n\)-dimensional manifolds, he carefully explained to his readers that he was introducing a new way of thinking, and we conclude with this quote from the introduction to his book (\cite{weyl1913}, p. V):
\begin{quote}
Eine solche
strenge Darstellung, die namentlich auch bei Begr\"{u}ndung der fundamentalen,
in die Funktionentheorie hineinspielenden Begriffe und S\"{a}reconntze der
%Analysis situs sich nicht auf anschauliche Plausibilit\"{a}t beruft, sondern
mengentheoretisch exakte Beweise gibt, liegt bis jetzt nicht vor. Die
wissenschaftliche Arbeit, die hier zu erledigen blieb, mag vielleicht als
Leistung nicht sonderlich hoch bewertet werden. Immerhin glaube ich
behaupten zu k\"{o}nnen, dass ich mit Ernst und Gewissenhaftigkeit nach
den einfachsten und Sachgem\"{a}ssesten Methoden gesucht habe, die zu dem
vorgegebenen Ziele f\"{u}hren; und an manchen Stellen habe ich dabei andere
Wege einschlagen m\"{u}ssen als diejenigen, die in der Literatur seit dem
Erscheinen von C. Neumanns klassischem Buche \"{u}ber „Riemanns Theorie
der Abelschen Integrale" (1865) traditionell geworden sind.\footnote{Such a rigorous presentation, which, namely by the establishing of the fundamental concepts and theorems in function theory and using theorems of the analysis situs which don't just depend on intuitive plausibility, but have set-theoretic exact proofs, does not exist. The scientific work that remains to be done   in this regard, may perhaps not be be particularly highly valued.  But, nevertheless, I believe I can maintain that I have tried in a serious and conscientious manner to find the simplest and most appropriate methods, that lead to the  asserted goal; and at many points I have had to proceed in a different manner than that which has become traditional in the literature since the appearance of C. Neumann's classical book about "Riemann's theory of Abeiian Integrals."}
\end{quote}

\bibliography{references}
\bibliographystyle{plain}%

\end{document}